\documentclass{amsart}
\usepackage{amssymb,latexsym,amscd}
\title[Quasiharmonic polynomials]{Quasiharmonic polynomials for Coxeter
groups and representations of Cherednik algebras}
\author[A.\,Berenstein]{Arkady Berenstein}
\address{Department of Mathematics, University of Oregon,
Eugene, OR 97403, USA}
\email{arkadiy@math.uoregon.edu}
\author[Yu.\,Burman]{Yurii Burman}
\address{Independent University of Moscow, 121002, 11, B.Vlassievsky per.,
Moscow, Russia}
\email{burman@mccme.ru}
\thanks{Research supported in part by the NSF (DMS) grants \#0102382 and 
\#0501103 (A.B.) and the RFBR grants \#N.Sh.4719.2006.1 and \#05-01-01012a 
(Y.B.)}

\date{}

\makeatletter

\RequirePackage{amssymb}

\newcommand{\theoremName}{Theorem}

\newcommand{\lemmaName}{Lemma}
\newcommand{\corollaryName}{Corollary}
\newcommand{\propositionName}{Proposition}
\newcommand{\conjectureName}{Conjecture}
\newcommand{\remarkName}{Remark}
\newcommand{\exampleName}{Example}
\newcommand{\definitionName}{Definition}
\newcommand{\proofName}{Proof}
\renewcommand{\proofname}{\proofName}

\newtheorem{theorem}{\theoremName}[section]

\newtheorem{proposition}[theorem]{\propositionName}

\newtheorem{conjecture}[theorem]{\conjectureName}
\newtheorem{corollary}[theorem]{\corollaryName}

\newtheorem{lemma}[theorem]{\lemmaName}

\theoremstyle{definition}

\newtheorem{definition}[theorem]{\definitionName}

\theoremstyle{remark}
\newtheorem{remark}[theorem]{\remarkName}

\newtheorem{example}[theorem]{\exampleName}

\let\@newpf\proof \let\proof\relax 
\def \namepf[#1] {\@newpf[\proofname\ #1]}
\newenvironment{proof}{\@ifnextchar[{\namepf}{\@newpf[\proofname]}}{\qed\endtrivlist}

\def \Real {{\mathbb R}}
\def \Integer {{\mathbb Z}}

\def \Complex {{\mathbb C}}

\def \lnorm#1\rnorm {\vphantom{#1}\left\|\smash{#1}\right\|}
\def \lmod#1\rmod {\vphantom{#1}\left|\smash{#1}\right|}

\newcommand \bydef {\stackrel{\mbox{\scriptsize def}}{=}}

\newcommand{\pder}[2] {\frac{\partial #1}{\partial #2}}

\newcommand{\pdertwo}[3] {\frac{\partial^2 #1}{\partial #2 \partial #3}}

\@ifundefined{name}{\newcommand \name[1] {\mathop{\rm #1}\nolimits}}%
{\relax}

\newcommand \eps {\varepsilon}
\renewcommand \phi {\varphi}
\renewcommand \rho {\varrho}

\numberwithin{equation}{section}

\let \CC=\Complex
\def \PP {\mathbb P}
\def \HH {\mathcal{H}}
\def \QH {\mathcal{QH}}
\def \Refl {\mathcal S}
\def \triv {\mathbf{1}}
\def \divides {\mathrel{\vert}}

\def \Sing {C^{\mathrm sing}}
\def \hilb {\name{hilb}}
\def \ch {\name{ch}}
\def \ev {\name{ev}}

\def \c {\sp{(c)}}

\begin{document}
 \begin{abstract}
We introduce and study deformations of finite-dimensional modules over
rational Cherednik algebras. Our main tool is a generalization of usual
harmonic polynomials for each Coxeter groups --- the so-called
quasiharmonic polynomials. A surprising application of this approach is
the construction of canonical elementary symmetric polynomials and their
deformations for all Coxeter groups.
 \end{abstract}

\maketitle

\tableofcontents

\section*{Introduction}

In this paper we introduce and study deformations of finite-dimensional
modules over the rational Cherednik algebras $H_c(W)$. Here $W$ is a finite
reflection group of the space $V$, and $c$ is a conjugation-invariant
complex-valued function on the set of reflections in $W$. The algebra
$H_c(W)$, first introduced in 1995 by I.\,Cherednik \cite{Cherednik}, is
a degenerate case of the double affine Hecke algebra and can be thought of
as a deformation of the cross product of the group algebra $\Complex W$ and
the symmetric algebra of $V^* \oplus V$; $c$ is the parameter of
deformation.

Rational Cherednik algebras attracted attention of many authors in the
recent decade. Finite-dimensional representations were studied in
\cite{BEG}; see also \cite{Chmutova}, \cite{CE}, \cite{Dez} and a review
\cite{Rouq}. Recently I.\,Gordon (\cite{GorQuot}, \cite{GorI},
\cite{GorII}) found remarkable connections between the theory of Cherednik
algebras and algebras of diagonal harmonics; the latter were introduced by
M.\,Haiman in \cite{Haiman} and have been the subject of intensive research
since then.

The paper is a part of a project having its goal to construct canonical
bases in the representations of the rational Cherednik algebra. In its
turn, this would allow the construction of a canonical basis in the algebra
of diagonal harmonics.

For generic $c$ the algebra $H_c(W)$ has no nontrivial finite-dimensional
modules. Certain  functions $c$ such that $|c| = r\in \Integer_{>0}$ (where
$|c|$  is, up to a multiple, the sum of all values of $c$ --- see
\eqref{Eq:abs c} for details) are exceptions to this rule: for such $c$ the
algebra $H_c(W)$ has a distinguished finite-dimensional module $A_r$. In
particular, this happens for $H_c(S_n)$ when $c = const. = r/n$ and $r$ and
$n$ are coprime (\cite{BEG}, see also \cite{CE} for a simpler proof), for
$H_c(I_2(m))$ when $m$ is odd, $c = const. = r/m\notin \Integer$, and when
$m$ is even, $c_1 + c_2 = r/m\notin \Integer$ (\cite{Chmutova}); for other
types of Coxeter groups see \cite{BEG}. The module $A_r$ for all these
cases has the form $M(\triv)/I_r$ where  $M(\triv)$ is the Verma-like
$H_c(W)$-module isomorphic to $S(V)$ as a vector space (first introduced by
C.\,Dunkl in \cite{DuInter}; see Section \ref{SSec:Defs} for exact
definition) and $I_r \subset M_r(\triv)$ is a proper maximal sub-module of
the form $I_r = V^{(r)} \cdot S(V)$. The subspace $V^{(r)} \subset S^r(V)$
(the space of {\em singular vectors}) is the common kernel of all {\em
Dunkl operators}.

In the present paper we propose  a certain flat $c$-deformation $V^{(r;c)}$
of the $W$-module $V^{(r)}$ such that $V^{(r;c)}$ is still a $W$-module and
the quotient algebra $A_{r;c} = S(V)/(V^{(r;c)}\cdot S(V))$ is
finite-dimensional and (yet conjecturally) specializes into $A_r$ when $c$
is exceptional as above. The advantage of the
deformed algebra $A_{r;c}$ is that it possesses an additional structure of
a flat family which, similarly to quantum groups, can help to understand
canonical bases for $A_r$.

A principal ingredient of the construction of $V^{(r;c)}$ is the concept of
{\em quasiharmonic elements}, i.e., elements of $S(V)$ killed by almost all
$W$-invariant combinations of Dunkl operators. Namely, $V^{(r;c)}$ is a
unique (except for the case when $W=D_{2m}$ - see Section \ref{SSec:ReprQH}
for details) quasiharmonic $W$-module in $S^r(V)$ isomorphic to $V^{(r)}$
(and thus we believe that $A_{r;c}$ is a canonical deformation of $A_r$). A
surprising application of this technique is a construction of a remarkable
system of canonical invariants $e_k^{(c)} \in S(V)^W$ of any Coxeter group
$W$. For $W = S_n$ these invariants are deformations of the elementary
symmetric polynomials $e_k$; a more detailed study of these invariants will
be the subject of a forthcoming paper.

The paper is organized as follows. In Section \ref{Sec:Deform} we give
necessary definitions (Section \ref{SSec:Defs}), describe the
quasiharmonic elements (Section \ref{SSec:Basic}), define the deformation
$V^{(r;c)}$, construct the invariants $e_k^{(c)}$ (Section
\ref{SSec:ReprQH}) and formulate several conjectures about properties of
$V^{(r;c)}$ and other related objects (Section \ref{SSec:Deform}). In
Section \ref{Sec:Dih} we verify these conjectures for the case $W = I_2(m)$.
Section \ref{SSec:Mult} deals with a more special question --- the algebra
$A_{r;c}$ and its Frobenius characteristic polynomial. Section
\ref{Sec:Apolar} is an appendix summarizing information about standard 
Frobenius algebras used elsewhere in the paper and containing some remarks 
we believe necessary for the future development of the subject.

\subsection*{Acknowledgments}
The authors are grateful  Michael Feigin, Alexander 
Veselov, Yuri Bazlov, and Sergey Yuzvinsky for stimulating discussions. Special thanks are due to Pavel Etingof for 
suggesting the idea of the proof of Theorem \ref{Th:SymmE}.  
Yu.\,Burman wishes to thank the University of Oregon, where this work was 
started, for its warm hospitality.

\section{Quasiharmonic elements: general case}\label{Sec:Deform}

\subsection{Definitions and notation}\label{SSec:Defs}

For a complex vector space $V$ we denote by $S(V)$ its symmetric algebra,
and by $\Complex[V]$ the algebra of polynomials functions on it. Clearly,
$\Complex[V] = S(V^*)$. Both $S(V)$ and $\Complex[V]$ are graded algebras;
their homogeneous components of degree $k$ are denoted $S^k(V)$ and
$\Complex_k[V]$, respectively.

Let $W$ be a finite reflection group of a Euclidean vector space $V_\Real$
such that the $W$-module $V_\Real$ is irreducible. Denote by $V \bydef
\Complex \otimes V_\Real$ the complexification of $V_\Real$. Clearly, $V$
is also an irreducible $W$-module; we will call it the {\em defining
$W$-module}. The complex space $V$ inherits from $V_\Real$ a symmetric
non-degenerate $W$-invariant form $(\cdot,\cdot)$. Denote by  $\Refl
\subset W$ the set of all reflections in $W$. To each $s \in \Refl$ one
associates a positive root $\alpha_s \in V$ and a coroot $\alpha_s^{{\vee}}
= 2(\alpha_s, \cdot)/(\alpha_s, \alpha_s) \in V^*$.

Fix a function $c: \Refl \to \Complex$ invariant under conjugation: $c(\tau
s \tau^{-1}) = c(s)$, and define the {\em rational Cherednik algebra}
$H_c(W)$ as follows. As a vector space, $H_c(W)$ is isomorphic to $S(V)
\otimes \Complex W \otimes S(V^*)$ where $\Complex W$ is the group algebra.
The multiplication is defined by the following requirements:

 \begin{enumerate}
\item The natural inclusions $S(V) \hookrightarrow H_c(W)$, $S(V^*)
\hookrightarrow H_c(W)$ and $\Complex W \hookrightarrow H_c(W)$ are algebra
monomorphisms. In particular, $x_1 x_2 = x_2 x_1$ for all $x_1, x_2 \in
S(V) \subset H_c(W)$, and the same is true for all $y_1, y_2 \in S(V^*)$.

\item The conjugation action of $\Complex W \subset H_c(W)$ on $V \subset
S(V) \subset H_c(W)$ and on $V^* \subset S(V^*) \subset H_c(W)$ is
isomorphic to the action of $\Complex W$ in the defining module $V$
and on its dual $V^*$, respectively: $wzw^{-1} = w(z)$ for any $w \in W$
and any $z \in V$ or $z \in V^*$.

\item For any $x \in V$, $y \in V^*$ one has
 \begin{equation}\label{Eq:Commut}
yx-xy = \langle y,x\rangle - \sum_{s \in \Refl} c(s) \langle
y,\alpha_s\rangle \langle \alpha_s^{{\vee}},x\rangle s.
 \end{equation}
Here brackets $\langle \cdot, \cdot\rangle$ mean pairing of $V^*$ and $V$;
the right-hand side lies in $\Complex W \subset H_c(W)$;
 \end{enumerate}

Evidently, the definition of $H_c(W)$ is sound, i.e.\ it does not depend on
the choice of a positive root $\alpha_s$ for every $s \in \Refl$.

We can collect all the algebras $H_c(W)$ together. Consider the
decomposition $\Refl = \bigsqcup_{i=1}^k C_i$ of the set of reflections
into conjugacy classes (in fact, $k=1$ or $2$ since $W$ is irreducible) and
denote $c_i \bydef c(s)$ for any $s \in C_i$. Then define a universal
algebra $H(W)$ over the ring $\Complex[c] = \Complex[c_1, \dots, c_k]$
exactly as in the previous paragraph with the field $\Complex$ replaced by
$\Complex[c]$ everywhere (including tensor multiplication). For each $c =
(c_1, \dots, c_k) \in \Complex^k$ we have a natural surjective
evaluation homomorphism of $\Complex$-algebras $\ev_c: H(W) \to H_c(W)$.
The algebra $H(W)$ is a convenient language to study the dependence of
various objects inside $H_c(W)$ on the parameter $c$. 

Define the projection $\pi_{12}: H_c(W) \to S(V) \otimes \Complex W$ by 
$\pi_{12}(q \otimes w \otimes p) = p(0) \cdot q \otimes w$ (where $p\mapsto 
p(0)$ is the evaluation of a polynomial $p$ at $0$). In particular, 
$\pi_{12}(q \otimes w \otimes p) = 0$ for every $p \in S(V^*)_+ = 
\bigoplus_{m > 0} S^m(V^*)$. For every $p \in S(V^*)$ define the {\em Dunkl 
operator} $\nabla_p\c: S(V) \otimes \Complex W \to S(V) \otimes \Complex W$ 
by the formula
 \begin{equation*}
\nabla_p\c(r) = \pi_{12}(pr).
 \end{equation*}
(in the right-hand side the product of elements $p \in S(V^*) \subset
H_c(W)$ and $r \in S(V) \otimes \Complex W \subset H_c(W)$ in the algebra
$H_c(W)$ is assumed).

 \begin{proposition}[cf.\ \cite{DO}]\label{Pr:DO}
The space $S(V) \otimes \Complex W$ possesses the following $H_c(W)$-module
structure: elements $q \in S(V)$ act there by multiplication in the first
factor, elements $w \in \Complex W$, as in the module $S(V) \otimes
\Complex W$, and elements $p \in S(V^*)$, by Dunkl operators. Explicitly,
the action of the Dunkl operator corresponding to a vector $y \in V^*$ on
an element $q \otimes w \in S(V) \otimes \Complex W$ is given by the
following formula:
 \begin{equation}\label{Eq:Dunkl}
\nabla_y\c (q \otimes w) = \pder{q}{y} \otimes w - \sum_{s \in \Refl}
c(s) \langle y,\alpha_s\rangle \frac{q - s(q)}{\alpha_s} \otimes s w.
 \end{equation}
 \end{proposition}

In \eqref{Eq:Dunkl} the symbol $\pder{q}{y}$ means the directional
derivative of a polynomial $q \in S(V) = \Complex[V^*]$. An element $s \in
\Refl \subset W$ acts on $q$ as in the module $S(V)$; the difference $q -
s(q) \in S(V)$ is divisible by $\alpha_s \in V$, so the right-hand side is
well-defined.

 \begin{proof}
Extend the notation $\nabla$ defining $\nabla_u\c(r) \bydef \pi_{12}(ur)$
for every $u \in H_c(W)$ and $r \in S(V) \otimes \Complex W$. Then show
that the correspondence $(u,r) \mapsto \nabla_u\c(r)$ is an $H_c(W)$-action
on $S(V) \otimes \Complex W$.

Indeed, for any $x \in H_c(W)$ one has $x - \pi_{12}(x) = \sum_i q_i
\otimes w_i \otimes p_i$ where $q_i \in S(V)$, $w_i \in \Complex W$ and
$p_i \in S^{k_i}(V^*)$ with $k_i > 0$. A similar decomposition takes place
for the element $u(x - \pi_{12}(x))$ where $u \in H_c(W)$, and therefore
$\pi_{12}(u(x - \pi_{12}(x))) = 0$ for all $u,x \in H_c(W)$. This implies
the equality $\nabla_{u_1u_2}\c(r) = \pi_{12}(u_1u_2r) = \pi_{12}(u_1
\pi_{12}(u_2r)) = \nabla_{u_1}\c(\nabla_{u_2}\c(r))$, as required.

Furthermore, the defining relation (\ref{Eq:Commut}) implies, by direct
computation, the following relation in $H_c(W)$:
 \begin{equation*}
(yq-qy)w = \pder{q}{y}\cdot w - \sum_{s \in \Refl} c(s) \langle
y,\alpha_s\rangle \frac{q - s(q)}{\alpha_s} \cdot  sw
 \end{equation*}
for all $y\in V^*$, $q\in S(V)$, $w\in \Complex W$. By definition of the
projection $\pi_{12}$, one has also $\pi_{12}(yqw) = \pi_{12}((yq-qy)w) =
(yq-qy)w$. Taking into account the factorization $S(V) \cdot \Complex W =
S(V) \otimes \Complex W$ (in $H_c(W)$), we obtain \eqref{Eq:Dunkl}.
 \end{proof}

One can replace $\Complex W$ in Proposition \ref{Pr:DO} by an arbitrary
$\Complex W$-module $\tau$. Following \cite{BEG} we denote this
$H_c(W)$-module by $M(\tau)$; it is isomorphic to $S(V)\otimes \tau$ as a
vector space. Modules $M(\tau)$ inherit the grading of $S(V)$: $M(\tau) =
\bigoplus_n M_n(\tau)$ where $M_n(\tau) = S^n(V) \otimes \tau$. In
particular, setting $\tau = \triv$ (the trivial $1$-dimensional
representation) one obtains a structure of the $H_c(W)$-module in $M(\triv)
= S(V)$ where elements of $S(V)$ act by multiplication, elements $w \in
\Complex W$, as in the module $S(V)$, and elements $p \in V^*$, by the
Dunkl operators in the original sense of \cite{DO}:
 \begin{equation}\label{Eq:DunklOrig}
\nabla_y\c(q) = \pder{q}{y} - \sum_{s \in \Refl} c(s) \langle
y,\alpha_s\rangle \frac{q - s(q)}{\alpha_s}.
 \end{equation}
Below we will refer to the operators \eqref{Eq:DunklOrig} as Dunkl
operators unless noted otherwise.

We can also consider a ``universal'' version of the above construction: if 
$\widetilde\pi_{12}: H(W) \to \Complex[c] \otimes S(V) \otimes \Complex W$ 
is the natural projection and $p \in S(V^*)$, then the universal Dunkl 
operator is defined by the formula $\nabla_p(r) \bydef 
\widetilde\pi_{12}(pr)$. It is clear that the ``universal'' and 
``specialized'' operators are related by the evaluation morphism: 
$\nabla_p\c\ev_c = \ev_c \nabla_p$; by a slight abuse of notation we will 
sometimes be writing $\nabla_p$ instead of $\nabla_p\c$.

Define a linear map $\eps: H_c(W) \to \Complex$ by $\eps(q \otimes w 
\otimes p) = p(0)q(0)$, where $p \mapsto p(0)$ and $q \mapsto q(0)$ are the 
evaluations of polynomials at $0$. In particular, $\eps(S(V^*)_+ \otimes w 
\otimes S(V)) = 0$ and $\eps(S(V^*) \otimes w \otimes S(V)_+) = 0$ for each 
$w\in W$. Define a bilinear pairing between the spaces $S(V^*)$ and $S(V)$ 
by the formula
 \begin{equation}\label{Eq:Pairing}
\langle p,q\rangle_c  = \eps(pq)
 \end{equation}
for $p \in S(V^*)$, $q \in S(V)$; so, $\langle p,q\rangle_c= \delta_{mn} 
\nabla_p\c(q)$ for $p\in S^m(V^*)$ and $q\in S^n(V)$. For $c = 0$ equation 
\eqref{Eq:Pairing} defines the usual pairing between $S(V^*)$ and $S(V)$.

For each function $c:\Refl\to \Complex$ denote
 \begin{equation}\label{Eq:abs c}
\lmod c\rmod = \frac{2}{\ell}\sum_{s \in \Refl} c(s).
 \end{equation}
where $\ell \bydef \dim V$. The following result was essentially proved in 
\cite[Proposition 2.1]{BEG}.

 \begin{proposition}
If $q \in V$ and $p \in V^*$, then $\langle p,q\rangle_c = (1-\lmod
c\rmod)\langle p,q\rangle_0$.
 \end{proposition}

The non-degenerate bilinear form $(\cdot, \cdot)$ on the space $V$ (used in
the definition of the reflection group $W$) defines an isomorphism of
vector spaces $\Phi: V \to V^*$. One can extend this to an isomorphism
$\Phi: S(V) \to S(V^*)$ and further, to the linear involution  $\phi: H_c(W) \to H_c(W)$ by $\phi(q \otimes w
\otimes p) = \Phi^{-1}(p) \otimes w^{-1} \otimes \Phi(q)$.

 \begin{proposition}[cf.\ \cite{BEG}]
The following is true:
 \begin{enumerate}
\item\label{It:AntiHom} The involution $\phi$ is an anti-automorphism of 
$H_c(W)$: $\phi(ab) = \phi(b) \phi(a)$.

\item\label{It:Symm} The pairing $\langle\cdot, \cdot\rangle_c$ is
symmetric with respect to $\phi$: $\langle p,q\rangle_c =
\langle\phi(q),\phi(p)\rangle_c$.

\item\label{It:Conj} The Dunkl operators are conjugate to multiplication
operators in $S(V^*)$: $\langle rp,q\rangle_c = \langle
r,\nabla_p(q)\rangle_c$ for any  $p,r \in S^k(V^*)$ and $q \in S^n(V)$.
 \end{enumerate}
 \end{proposition}

 \begin{proof} 
Since $H_c(W)$ is spanned by elements of the form $qwp$ where $q \in S(V)$, 
$w \in W$ and $p \in S(V^*)$, in order to prove assertion 
\eqref{It:AntiHom} it suffices to verify that
 \begin{equation}\label{Eq:3case}
\phi(aqwp) = \Phi^{-1}(p)w^{-1}\Phi(q) \phi(a)
 \end{equation}
for $a \in V\cup  W \cup V^*$. For $a \in V$ one has $\phi(a) = \Phi(a)$, and \eqref{Eq:3case} is true 
because $\Phi: S(V) \to S(V^*)$ is an algebra homomorphism. For $a \in W$ one has $\phi(a) = a^{-1}$, and therefore
$\phi(aqwp) = \phi(a(q)awp) = \Phi^{-1}(p) w^{-1} a^{-1} \Phi(a(q))$.
The group $W$ acts in $S(V)$ by $(\cdot,\cdot)$-orthogonal operators, so
that $\Phi(a(q)) = a(\Phi(q))$ and therefore
 \begin{equation*}
\phi(aqwp) = \Phi^{-1}(p) w^{-1} a^{-1} a(\Phi(q)) = \Phi^{-1}(p) w^{-1}
\Phi(q) a^{-1}.
 \end{equation*}

Let now $a \in V^*$. It follows from \eqref{Eq:Commut} that for every $b
\in V$ one has
 \begin{align*}
\phi([a,b]) &= \phi\bigl(\langle a,b\rangle - \sum_{s \in \Refl} c(s)
\langle a,\alpha_s\rangle \langle \alpha_s^{{\vee}},b\rangle s\bigr) =
\langle a,b\rangle - \sum_{s \in \Refl} c(s) \langle a,\alpha_s\rangle
\langle \alpha_s^{{\vee}},b\rangle s \\
&= \langle \Phi(b),\Phi^{-1}(a)\rangle - \sum_{s \in \Refl} c(s) \langle 
\Phi(b),\alpha_s\rangle \langle \alpha_s^{{\vee}},\Phi^{-1}(a)\rangle s = 
[\Phi(b),\Phi^{-1}(a)].
 \end{align*}
Trivial induction shows that $\phi([a,q]) = [\phi(q),\phi(a)]$ for every $q
\in S(V)$.

Now one has $aqwp = qawp + [a,q]wp = qw\cdot w^{-1}(a)p + [a,q]wp$. Since
 \begin{equation*}
\phi([a,q]wp) = \Phi^{-1}(p) w^{-1} \phi([a,q]) = \Phi^{-1}(p) w^{-1}
[\phi(q),\phi(a)]
 \end{equation*}
and also,
 \begin{align*}
\phi(qw\cdot w^{-1}(a)p) &= \Phi^{-1}(w(a)p) w^{-1} \Phi(q)\\
&= \Phi^{-1}(p) w^{-1}(\Phi^{-1}(a)) w^{-1} \Phi(q) = \Phi^{-1}(p) w^{-1} 
\phi(a) \phi(q),
 \end{align*}
one has eventually
 \begin{equation*}
\phi(aqwp) = \Phi^{-1}(p)w^{-1}(\phi(a) \phi(q) + [\phi(q),\phi(a)]) =
\phi(qwp) \phi(a),
 \end{equation*}
so that assertion \eqref{It:AntiHom} is proved.

To prove assertion \eqref{It:Symm} note that the involution $\phi$ is 
graded, hence $\eps(x - \phi(x)) = 0$ for all $x \in H_c(W)$. Therefore, 
assertion \eqref{It:AntiHom} implies that
 \begin{equation*}
\langle \phi(q), \phi(p)\rangle_c = \eps(\phi(q)\phi(p)) = \eps(\phi(pq)) =
\eps(pq) = \langle p,q\rangle_c.
 \end{equation*}

Assertion \eqref{It:Conj} follows from Proposition \ref{Pr:DO}:
 \begin{equation*}
\langle rp,q\rangle_c = \nabla_{rp}\c(q) = \nabla_r\c(\nabla_p\c(q)) = 
\langle r,\nabla_p\c(q)\rangle_c
 \end{equation*}
This finishes the proof of the proposition.
 \end{proof}

 \begin{definition}
An element $x\in S(V)$ is {\em $(W,c)$-harmonic} if $\nabla_p\c(x) = 0$ for
any $W$-invariant polynomial $p\in S(V^*)^W_+$. An element $x \in 
\Complex[c] \otimes S(V)$ is {\em $W$-harmonic} if $\nabla_p(x) = 0$ for 
any $W$-invariant polynomial $p\in S(V^*)^W_+$.
 \end{definition}

A classical Chevalley theorem (see \cite{Chev}) says that $S(V^*)^W$ is a
polynomial algebra of $\ell$ variables (recall that $\ell = \name{rk}(W) =
\dim V$). These variables correspond to some invariant polynomials ({\em
elementary invariants}). The elementary invariants can be chosen in many
ways (see, though, Proposition \ref{Pr:unique elementary invariant} but their
degrees are uniquely defined and called {\em exponents} of the group $W$.
Denote the exponents as $d_1 \le \dots \le d_\ell = h$ (the largest
exponent $h$ is the Coxeter number); and let $e_{d_1}, \dots, e_{d_\ell} =
e_h$ be some elementary invariants of $W$ in $S(V^*)$ of degrees
$d_1,\dots,d_\ell$ respectively.

Using this description of invariants one can give a simpler
characterization of $(W,c)$-harmonic elements:

 \begin{proposition}
The following are equivalent for $x \in S(V)$ (respectively, $x \in
\Complex[c] \otimes S(V)$):
 \begin{enumerate}
\item $x$ is $(W,c)$-harmonic (respectively, $W$-harmonic);

\item $\nabla_p(x) = 0$ for any $W$-invariant $p \in S(V^*)^W_+$ such
that $\deg(p)\le h$.

\item $\nabla_{e_{d_i}}(x)=0$ for $i=1,2,\dots,\ell$.
 \end{enumerate}
 \end{proposition}

The proof is obvious.

Since $W$ is a {\em real} reflection group (a Coxeter group), one always 
has the inequality $h = d_\ell > d_{\ell-1}$ (\cite{Hum}), i.e.\ the 
invariant $e_h$ is unique modulo invariants of smaller degrees. By removing 
$e_h$ from the consideration, we obtain the following object:

 \begin{definition}
An element $x\in S(V)$ is $(W,c)$-{\em quasiharmonic} if $\nabla_p\c(x)=0$
for any $W$-invariant polynomial $p\in S(V^*)^W$ such that $\deg(p) < h$.
An element $x\in \Complex[c] \otimes S(V)$ is $W$-{\em quasiharmonic} if
$\nabla_p(x)=0$ for any such polynomial. Or, equivalently, $x$ is
$(W,c)$-quasiharmonic ($W$-quasiharmonic) if $\nabla_{e_{d_i}}(x) = 0$ for
$i = 1,2,\dots,\ell-1$.
 \end{definition}

Denote by $\HH\c$ the space of $(W,c)$-harmonic elements and by $\QH\c
\supset \HH\c$ the space of $(W,c)$-quasiharmonic elements. These spaces
are graded:
 \begin{equation*}
\HH\c = \bigoplus_{n \ge 0} \HH_n\c, \quad \QH\c =
\bigoplus_{n \ge 0} \QH_n\c
 \end{equation*}
where $\HH_n\c = \HH\c \cap S^n(V)$ and $\QH_n\c = \QH\c \cap S^n(V)$. The 
same is true for the ``universal'' versions of these spaces, $\HH \subset 
\QH$. Clearly, $\ev_c\HH \subseteq \HH\c$ and $\ev_c\QH \subseteq \QH\c$.

 \begin{remark}\label{rem:complex reflections}
The notions and results of this section are also applicable to {\em 
complex} reflection groups. The reason for this is two-fold. First, the 
algebra of invariants of such a group is free (see \cite{SheTodd}), and the 
exponents $d_1 \le d_2 \le \dots \le d_\ell$ are well-defined. Second, the 
inequality $d_\ell > d_{\ell-1}$ holds for almost all groups. More 
precisely, it fails only for the following irreducible complex reflection 
groups of rank $\ell = 2$ with $(d_1,d_2) = (d,d)$ (see e.g., 
\cite{BrMaRo}):

 \begin{enumerate}
\item The series $G(m,2,2)$, $m>1$, with $d=m$.

\item The exceptional groups $G_7$, $G_{11}$, $G_{19}$ with $d=12, 24, 60$, 
respectively. 
 \end{enumerate}
 \end{remark}

\subsection{Basic properties of quasiharmonic elements}\label{SSec:Basic}

Following \cite{DunklJeuOpdam}, denote by $\Sing$ the set of all 
$W$-invariant functions $c: \Refl \to \Complex$ such that there exists an 
element $x \in S(V)$ of positive degree such that $\nabla_y\c(x) = 0$ for 
all $y \in V^*$. Elements of $\Sing$ are called {\em singular functions}; 
other conjugation-invariant functions on $\Refl$ are called {\em regular}. 
Elements $x \in S(V)$ annihilated by all the Dunkl operators are called 
{\em singular vectors}. It is easy to see that the group $W$ acts on the 
set $X$ of singular vectors, and that the space $S(V)X \subset S(V) = 
M(\triv)$ is a $H_c(W)$-submodule. Thus, $c \in \Sing$ if and only if the 
$H_c(W)$-module $M(\triv)$ is reducible (see \cite{BEG} for details). It 
was proved in \cite{DunklJeuOpdam} that the function $c = const.$ is 
singular if and only if $c = r/d_i$ where $d_i$ is an exponent of $W$ and 
$r$ is a positive integer not divisible by $d_i$. In particular, $c = 0$ is 
a regular function. Conjugation-invariant functions $c \ne const.$ (taking, 
actually, two different values) exist for Coxeter groups of types $B_n$, 
$F_4$, and $I_2(2m)$; see \cite{DunklJeuOpdam} for description of $\Sing$ 
in these cases.

 \begin{proposition}[\cite{DunklJeuOpdam}] \label{Pr:Intertw}
For all regular $c$ the $S(V^*)$-actions on $M(\triv) = S(V)$ induced by
the $H_c(W)$-action are isomorphic to each other. In particular, they are
isomorphic to the natural $S(V^*)$-action on $S(V)$ by differential
operators with constant coefficients (this corresponds to $c = 0$).
 \end{proposition}

 \begin{proposition} \label{Pr:deformed hilbert series}
For each regular $c: \Refl \to \Complex$ the Hilbert series of $\QH\c$
and of  $\HH\c$ are given by the following formulas:
 \begin{equation*}
\hilb(\QH\c,t) = \hilb(\HH\c,t)/(1-t^h)=\frac{(1-t^{d_1})\dots
(1-t^{d_{\ell-1}})}{(1-t)^\ell},
 \end{equation*}
 \end{proposition}

 \begin{proof}
Let $I \subset S(V^*)$ be a homogeneous ideal. Since the pairing $\langle
\cdot, \cdot\rangle_c$ is non-degenerate, there is an isomorphism of graded
vector spaces:
 \begin{equation*}
S(V^*)/I \cong I^\perp,
 \end{equation*}
where $I^\perp = \{x \in S(V) \mid \langle I,x\rangle_c = 0\}$ is the
orthogonal complement of $I$. This implies $\hilb(S(V^*)/I,t) =
\hilb(I^\perp,t)$.

Recall that a homogeneous ideal $I \subset S(V^*)$ is called {\em free} if
it admits a free Koszul resolution:
 \begin{equation}\label{Eq:Koszul}
0\to S(V^*)\otimes \Lambda^k(U)\to \dots \to S(V^*)\otimes \Lambda^2(U)\to
S(V^*)\otimes U\to S(V)\to S(V^*)/I \to 0
 \end{equation}
where $U \subset S(V^*)$ is a $k$-dimensional space of homogeneous
polynomials that generates the ideal $I$.

One has $\HH\c = I_1^\perp$ and $\QH\c = I_2^\perp$ where ideals
$I_1$ and $I_2$ are both free. Their generating sets $U_1$ and $U_2$ are
linear spans of $e_{d_1}, \dots, e_{d_k}$ where $k = \ell$ for $U_1$ and $k
= \ell-1$ for $U_2$. Then \eqref{Eq:Koszul} gives
 \begin{equation*}
\hilb(S(V^*)/I,t) = \frac{(1-t^{d_1})\cdots (1-t^{d_k})}{(1-t)^\ell}.
 \end{equation*}
and the result follows.
 \end{proof}

 \begin{corollary}
For any regular $c$ the following is true:
 \begin{enumerate}
\item $\QH\c_r = \HH\c_r$ for all $r < h$.

\item $\dim \QH\c_{h+1} = \dim \HH\c_{h+1} + \ell$.

\item For $r \gg 1$ the dimension of $\QH\c_r$ does not depend on $r$
and is equal to $\frac{\lmod W\rmod}{h} = d_1 d_2 \dots d_{\ell-1}$.
 \end{enumerate}
 \end{corollary}

\subsection{$W$-action on the space of quasiharmonics}\label{SSec:ReprQH}

The Proposition \ref{Pr:deformed hilbert series} above admits a
representation-theoretic refinement. For any graded $W$-submodule $U =
\bigoplus_{n \ge 0} U_n \subset S(V)$ and an irreducible $W$-module $\tau$
define the graded character $\ch_\tau(U,t)$ by:
 \begin{equation*}
\ch_\tau(U,t) = \sum_{n \ge 0} [\tau:U_n] t^n
 \end{equation*}
where $[\tau:U_n]$ is the multiplicity of $\tau$ in $U_n$.

All the constructions in the proof of Proposition \ref{Pr:deformed hilbert
series} are
$W$-equivariant, and each graded component of the ideals involved is a
$W$-module. So we immediately obtain 

 \begin{proposition}\label{Pr:Hilb} For any simple $W$-module $\tau$ and any regular $c$ one has:
 $$\displaystyle{\ch_{\tau}(\QH\c,t) = \frac{\ch_{\tau}(\HH\c,t)}{1-t^h}} \ .$$
  In particular, if
 $\tau = \triv$ is the trivial $1$-dimensional $W$-module, one has:
  \begin{equation}
  \label{Eq:ChTriv}
\ch_{\triv}(\QH\c,t) = 1/(1-t^h) = 1 + t^h + t^{2h} + \dots
 \end{equation}
  \end{proposition}

\subsubsection{Invariants}\label{SSSec:Invar} 
We obtain the following direct corollary of Proposition \ref{Pr:Hilb}. 
 \begin{corollary} \label{cor:multiple invariants}
For regular $c$ the degree of each homogeneous quasiharmonic invariant $q \in \QH^W$ is 
divisible by $h$. For every $k = 0,1, \dots$ there is exactly one, up to 
proportionality, quasiharmonic invariant $e_{kh}\c$ of degree $kh$.
 \end{corollary}
 
 \begin{proof} 
According to  \eqref{Eq:ChTriv}, for each regular $c$ and each 
$k\ge 0$, 
 \begin{equation*}
\dim \QH_d^W = \begin{cases}
1  & \text{if $h$ divides $d$,}\\
0  & \text{otherwise.}
 \end{cases}
 \end{equation*}
Therefore, for each $k>0$ there exists a unique (up to a multiple) 
$W$-invariant $u_{kh} \in \Complex(c)\otimes \QH_{kh}$. Multiplying this 
invariant by an appropriate polynomial in $\Complex[c]$, we obtain 
$e_{kh}$. 
 \end{proof}

Polynomials $e_{kh}\c$ for the case $W = I_2(m)$ (the dihedral group) are 
computed in Section \ref{SSec:QH}. See also Section \ref{SSec:Comput} and 
for some explicit formulas in the case $W = S_4$.

The Dunkl operators commute, and therefore they act on the space $\QH\c$ of
quasiharmonics making it a $\Complex W \ltimes S(V^*)$-module (where
$\Complex W \ltimes S(V^*)$ denotes the cross-product of $\Complex W$ and
$S(V^*)$ isomorphic to $\Complex W \otimes S(V^*)$ as a vector space).

 \begin{theorem} \label{Th:SymmE}
For regular $c$ the $\Complex W \ltimes S(V^*)$-module $\QH\c$ is generated
by the elements $e_{kh}\c$, $k = 0,1,2,\dots\,$. More precisely,
 \begin{equation}\label{Eq:harmonic generation}
\QH\c = \{\nabla_p\c(e_{kh}\c) \mid p \in \HH^*, k = 0,1,2,\dots\}.
 \end{equation}
where $\HH^* \subset S(V^*)$ is the space of harmonic elements of $S(V^*)$.
 \end{theorem}

To prove the theorem we will need some auxiliary results. Let $I^* 
\subset S(V^*)$ be the ideal generated by all $p\in S(V^*)_+^W$ of degree less than $h$. 

 \begin{lemma} \label{Lm:Radical}
The ideal $I^*$ is radical.
 \end{lemma}

 \begin{proof}
 One has a harmonic decomposition $S(V^*) = S(V^*)^W \otimes \HH$. Therefore, if we fix generators $e_{d_1}, \dots, e_{d_{\ell-1}}, e_{d_\ell} = e_h$ of $S(V^*)^W$ (so that  $S(V^*)^W= \Complex[e_{d_1}, \dots, 
e_{d_\ell}]$, we see that 
algebra $A \bydef S(V^*)/I^*$ is spanned by all the elements $\overline x \cdot \overline {e_h}^k$ for $x \in \HH$ and $k = 0,1,\dots$ (here $\overline x$ denotes the image of $x$ in $A$). On the other hand, Proposition 
\ref{Pr:deformed hilbert series} ensures that Hilbert series of $A$ and 
$\HH[e_h]=\HH\otimes \Complex[e_h]$ are the same, so that and $A \cong \overline \HH[\overline {e_h}]$ (as a vector space) and  $\overline \HH\cong \HH$. Hence, 
no polynomial of $\overline {e_h} \in A$ with constant coefficients is a zero divisor 
in $A$. 

The system of equations 
$e_{d_1} = \dots = e_{d_{\ell-1}} = 0, e_h = 1$ has $\lmod W\rmod = d_1 
\dots d_\ell$ distinct solutions. It means that the hypersurfaces defined respectively by $e_i = 0$, $i 
= 1, \dots, \ell-1$ and $e_h = 1$ intersect transversally in their common 
points (the derivatives are linearly independent). Therefore, the algebra $B \bydef A/(\overline{e_h} - 1) = 
S(V^*)/(e_{d_1}, \dots, e_{d_{\ell-1}}, e_h - 1)$ is semisimple of dimension $|W|$   hence contains no 
non-zero nilpotents.

Suppose that the ideal $I^*$ is not radical, that is, $A$ contains 
nilpotent elements; denote by $z \in A$ a nilpotent element of the smallest 
degree (as a polynomial in $\overline{e_h}$). Since the image of $z$ under the natural 
projection $A \to B$ is zero, one has $z = (e_h - 1) y$ for some  $y \in A$. Then $0 = z^N = (\overline{e_h} - 1)^N y^N$ for some $N$. 
Since $(\overline{e_h} - 1)^N$ is not a zero divisor, this implies $y^N = 0$ contrary 
to the assumption that the degree of $z$ is minimal. The lemma is proved.
 \end{proof}

Consider now the exponential map $\exp(x) \bydef \sum_{n \ge 0} 
\frac{x^n}{n!}$ from $V$ to the completion $\widehat{S(V)}$ of $S(V)$. 

 \begin{lemma} \label{Lm:PropExp}
 \begin{enumerate}
\item\label{It:DiffExp}  $\partial_p \exp(v) = p(v) \exp(v)$ for any $v \in V$ and $p \in S(V^*)$,   
where $\partial_p \bydef \nabla^{(0)}_p$ is the differential operator with constant coefficients on 
$\widehat{S(V)}$ corresponding to the polynomial $p$.

\item\label{It:LinIndepExp} If $v_1, \dots, v_N \in V$ are all distinct 
then $\exp(v_1), \dots, \exp(v_N) \in \widehat{S(V)}$ are linearly 
independent.
 \end{enumerate}
 \end{lemma}

 \begin{proof}
The first assertion is obvious. To prove the second one, consider a linear 
functional $\phi: V \to \Complex$ such that $\lambda_1 \bydef \phi(v_1), 
\dots, \lambda_N \bydef \phi(v_N)$ are all distinct. Extend $\phi$ 
naturally to $\widehat{S(V)}$, then $\phi(\exp(v)) = \exp(\phi(v))$ for all 
$v \in V$. If $\sum_{i=1}^N \alpha_i \exp(v_i) = 0$ is a linear dependence 
then one must have $\sum_{i=1}^N \alpha_i \exp(\lambda_i t) = 0$ for any $t 
\in \Complex$. But functions $\exp(t \lambda_i)$ with distinct $\lambda_i$ 
are obviously linearly independent.
 \end{proof}

 \begin{proposition}\label{Pr:annihilator}
Let
 \begin{equation}\label{Eq:tilde quasiharmonics}
Q \bydef \{\partial_p (e_{kh}^{(0)}) \mid p \in S(V^*), k = 0,1,\dots\}
 \end{equation} 
where $\partial_p$ is as in Lemma \ref{Lm:PropExp}. Then $Q^\perp = I^*$ 
where the annihilator is computed with respect to the standard pairing 
between $S(V)$ and $S(V^*)$.
 \end{proposition}

 \begin{proof}
For each $v \in V$ define an element $F(v) \in \widehat{S(V)}$ by the 
formula
 \begin{equation*}
F(v) = \sum_{w \in W} \exp(w(v))
 \end{equation*}
Then $\partial_p F(v) = \sum_{w \in W} p(w(v)) \exp(w(v))$ by assertion 
\ref{It:DiffExp} of Lemma \ref{Lm:PropExp}. Assertion \ref{It:LinIndepExp} 
of the same lemma implies that $\partial_p F(v) = 0$ if and only if 
$p(w(v)) = 0$ for all $w \in W$. Take $v_0$ such that $p(v_0) = 0$ for any 
$p \in I^*$; then one has
 \begin{equation*}
F(tv_0) = \sum_{k\ge 0} e_{kh}^{(0)} t^k
 \end{equation*}
for $t \in \Complex$ (under certain normalization of coefficients). By 
Lemma \ref{Lm:Radical}, $\partial_p F(tv_0) = 0$ if and only if $p \in 
I^*$. Hence $I^*$ is the annihilator of $Q$. 
 \end{proof}

 \begin{proof}[of Theorem \ref{Th:SymmE}]
By Proposition \ref{Pr:Intertw} it suffices to prove the theorem for only one
regular $c$. We will do this for $c = 0$. Graded $S(V^*)$-modules 
$\QH^{(0)}$ and $Q$ have the same annihilator $I^*$ in $S(V^*)$. Since $Q 
\subseteq \QH^{(0)}$ we immediately obtain $\QH^{(0)} = Q$. This proves the 
first assertion of the theorem.

Equation \eqref{Eq:harmonic generation} holds because $S(V^*) = S(V^*)^W 
\otimes \HH^*$ and, clearly, $\nabla_{e_h}\c (e_{kh}\c) \in \Complex \cdot 
e_{(k-1)h}\c$ while $\nabla_p\c (e_{kh}\c) = 0$ for $p\in S^r(V^*)^W$, $r < 
h$ by definition.

The theorem is proved.
 \end{proof}

Denote now $\QH^{(c,d)}$ the set of all elements $x \in S(V)$ such that 
$\langle p, x\rangle_c = 0$ for any $W$-invariant polynomials $p \in 
S^k(V^*)$ of degree $k < d$ (so that $\QH\c = \QH^{(c,h)}$). Clearly, 
$\QH^{(c,d)}$ is graded: $\QH^{(c,d)} = \bigoplus_{n \ge 0} \QH_n^{(c,d)}$. 
The ``universal'' versions of these spaces will be denoted 
$\QH^{(\cdot,d)}$ and $\QH^{(\cdot,d)}_n$, respectively.

Reasoning like in the proof of Proposition \ref{Pr:deformed hilbert series} and in Proposition \ref{Pr:Hilb},
one obtains 

 \begin{proposition} \label{Pr:HSTrunc}
For each simple $W$-module $\tau$ and any regular $c$ and $i = 1,\dots,\ell$:
 \begin{equation}\label{Eq:HSTrunc}
\hilb_\tau(\QH^{(c,d_i)},t) = \frac{\hilb_\tau(\HH\c,t)}{(1-t^{d_i}) \dots
(1-t^{d_{\ell}})}.
 \end{equation}
  In particular, if
 $\tau = \triv$ is the trivial $1$-dimensional $W$-module, one has:
 \begin{equation}\label{Eq:HSTruncTriv}
\hilb_{\bf 1}(\QH^{(c,d_i)},t) = \frac{1}{(1-t^{d_i}) \dots
(1-t^{d_{\ell}})} \ .
 \end{equation}
 \end{proposition}

This implies the following important assertion:

 \begin{proposition}\label{Pr:unique elementary invariant}
For each $i = 1,2,\dots,\ell$ and each regular $c$ one has
 \begin{equation}
 \label{Eq:Dim1Inv}
\dim (\QH^{(c,d_i)}_{d_i})^W = \begin{cases}
1, &d_i < d_{i+1},\\
2, &d_i = d_{i+1},
 \end{cases}
 \end{equation}
Therefore, if all $d_i$ are distinct, then for every $i = 1,2,\dots,\ell$ there
exists a $W$-invariant $e\c_{d_i} \in \QH^{(\cdot,d_i)}_{d_i}$ unique up
to proportionality.
 \end{proposition}
 
 \begin{proof} Clearly, \eqref{Eq:Dim1Inv}  a direct consequence of  \eqref{Eq:HSTruncTriv}. In its turn, if all $d_i$ are distinct, 
 \eqref{Eq:Dim1Inv} implies that for each
$i\le \ell$ there exists a unique (up to a multiple) $W$-invariant $u_{d_i} 
\in \Complex(c) \otimes \QH^{(\cdot,d_i)}_{d_i}$. Multiplying this 
invariant by an appropriate polynomial in $\Complex[c]$, we obtain 
$e\c_{d_i}$. The proposition is proved.
 \end{proof}

 \begin{remark} \label{Rm:DefEC}
Proportionality in this corollary means possibility to multiply $e\c_{d_i}$ 
by any polynomial in $\Complex[c]$. We can assume, though, that $e\c_{d_i} 
\in \Complex[c] \otimes S(V)$ is an irreducible polynomial; $e\c_{d_i}$ is 
then unique up to the multiplication by a constant. When it does not lead 
to confusion, we will denote by the same symbol $e^{(c_0)}_{d_i}$ the 
evaluation $\ev_{c_0}(e\c_{d_i})$ for a particular function $c_0:\Refl\to 
\Complex$.
 \end{remark}

 \begin{remark}
The case $d_i = d_{i+1}$ for irreducible Coxeter groups occurs only if $W = 
D_\ell$ with even $\ell$ and $i = \ell/2$.  In this case there exists an 
involutive automorphism $\sigma$ of $V$ acting non-trivially on the 
invariants in $\QH^{(\cdot,d_i)}_{d_i}$. We denote $e\c_{d_i}$ the 
invariant of $\sigma$, and $e\c_{d_i,-}$, its skew-invariant. 

In complex reflection groups the phenomenon $d_i = d_{i+1}$ is also very 
rare: in addition to those groups listed in Remark \ref{rem:complex 
reflections}, it occurs only for $G(kp,p,\ell)$ with $i = \ell/p$ provided 
$p$ divides $\ell$. Therefore, the invariants $e\c_{d_i}$ make sense for 
almost all complex reflection groups, too.
 \end{remark}

The polynomials $e\c_{d_i}$ should be interesting to study. In particular, 
$e_h\c$ is quasiharmonic and $\QH\c_h = \HH\c_h \oplus \Complex e_h\c$.

 \begin{theorem} \label{Th:BetaIndep}
For generic $c$ one has $S(V)^W = \Complex[e\c_{d_1}, \dots,
e\c_{d_\ell}]$. In other words, $e\c_{d_i}$, $i = 1, \dots, \ell$ form a
system of elementary $W$-invariants.
 \end{theorem}

 \begin{proof}
Let first $c = const.$, and let $P(c,e\c_{d_1}, \dots, e\c_{d_\ell}) = 0$
be an algebraic dependence between $e\c_{d_1}, \dots, e\c_{d_\ell}$.
Without loss of generality $P(0,x_1, \dots, x_n) \not\equiv 0$ --- else we
could have canceled the extra $c^k$ factor. But $e_{d_i}^{(0)}$ are
elementary invariants of the group $W$ in the sense the classical
Chevalley's theorem \cite{Chev}, so they are algebraically independent.

Let now $c \ne const.$; in this case (possible for Coxeter groups $B_n$, 
$I_2(m)$ with $m$ even, and $F_4$) $c$ assumes two values, $c_1$ and $c_2$. 
Let $P(c_1, c_2 ,e\c_{d_1}, \dots, e\c_{d_\ell}) = 0$ be an algebraic 
dependence between $e\c_{d_1}, \dots, e\c_{d_\ell}$. Without loss of 
generality $P(a,a,x_1, \dots, x_n) \not \equiv 0$ --- else we could have 
canceled the extra $(c_1-c_2)^k$ factor. But $\left. 
e_{d_i}^{(c)}\right|_{c_1=c_2=a}$ are algebraically independent for generic 
$a \in \Complex$, as we already proved above.
 \end{proof}

 \begin{remark}
Apparently, the first construction (recursive) of elementary invariants 
belongs to Dynkin (see e.g.\ \cite{Kos}) who used a similar approach for 
constructing the elementary $\mathrm{Ad}$-invariant polynomials over 
semisimple Lie algebras.
 \end{remark}

\subsubsection{Defining module}

Consider now Proposition \ref{Pr:Hilb} for $\tau = V$, the defining
$W$-module. It is well-known (see \cite{Hum}) that
 \begin{equation}\label{Eq:standard multiplicity}
\ch_V(\HH^{(0)},t) = t^{d_1-1} + t^{d_2-1} + \dots + t^{d_\ell-1}.
 \end{equation}
By Proposition \ref{Pr:Intertw}, the same formula is true for every regular
$c$, so that for $r \ge 1$ one has
 \begin{equation}\label{Eq:MutDef}
[V:\QH\c_r]=\begin{cases}
0 & \text{if $r \bmod h \notin \{d_1-1, d_2-1, \dots, d_\ell-1\}$},\\
2 & \text{if $r \bmod h = d_{\ell/2}-1 = d_{\ell/2 + 1}-1
= \frac{h}{2}-1$},\\
1 & \text{otherwise}.
 \end{cases}
 \end{equation}

Thus, for regular $c$ and for every $r \equiv d_i-1 \bmod h$, $i =
1,2,\dots,\ell$, in the case $d_i \ne d_{i \pm 1}$ the space $\QH_r\c$
contains a unique copy of the defining module $V$, which we denote
by $V^{(r;c)}$. In the case $r \equiv d_{\ell/2} - 1 = d_{\ell/2+1} - 1 =
h/2 - 1$ (it occurs only in $D_\ell$ type with even $\ell$) there is an
involutive automorphism $\sigma$ of $V$ acting non-trivially on
$\QH_r\c$. Then we take $V^{(r;c)}$ to be the $\sigma$-invariant copy
of $V$ in $\QH_r\c$.

More generally, let $\tau$ be an irreducible $W$-module such that
$[\tau:\QH_r^{(0)}] = 1$. By Proposition \ref{Pr:Intertw}, $[\tau:\QH_r\c] =
1$ for any regular $c$; we denote the unique copy of $\tau$ in $\QH_r\c$ by
$\tau^{(r;c)}$.

Theorem \ref{Th:SymmE} suggests the following procedure of computing 
$V^{(r;c)}$ where $r \equiv d_i-1 \mod h$ for some $i\le \ell$.  For each 
$j$ denote by now $V^*_j$ the unique copy of $V^*$ in the  space 
$\HH_{d_j-1}^*$ of the $W$-harmonic polynomials in $S(V^*)$.

 \begin{corollary}[of Theorem \ref{Th:SymmE}]
For all regular $c$ and any $r>0$ of the form $r = kh+d_i-1$ one has:
 \begin{equation*}
V^{(r;c)}=\nabla_{V^*_{\ell+1-i}}\c (e_{(k+1)h}\c)
 \end{equation*}
 \end{corollary}
\noindent (recall that $d_i+d_{\ell+1-i} = h+2$ because $W$ is  a real 
reflection group).

More generally, for a submodule $U \subset \HH^*_d$, $1 \le d < h$, one can 
consider a submodule $U^{(r;c)} \bydef \nabla_U(e_{r+d}\c) \subset 
\QH_r\c$, where $r = kh - d$. Theorem \ref{Th:SymmE} implies that 
$U^{(r;c)}$ is isomorphic to $U$.

\subsection{Flat deformations}\label{SSec:Deform}

We formulate now some conjectures about structure of  the
space of universal quasiharmonics $\QH$ and other related objects. The
next section (Section \ref{SSec:Comput}) contains various computations and
other evidence supporting the conjectures.

 \begin{conjecture}\label{Cj:SymmE}
For every singular function $c$ all the singular vectors of $M\c(\triv)$
belong to the specialization $\ev_c\QH$ of $\QH \subset H(W)$.
 \end{conjecture}

In other words, we conjecture that the space of singular vectors in
$M\c(\triv)$ is ``deformable'' in the class of quasiharmonics. We will
prove this conjecture for dihedral groups in Section \ref{SSec:QH} (for
singular $c = const.$ and generic singular $c \ne const.$). If $W = S_n$ is
a symmetric group, then explicit calculations for small $n$ in low degrees
of $M(\triv)$ support the conjecture; see Section \ref{SSec:Comput} below.

Suppose now that $r \equiv d_i-1 \bmod h$ for some $i = 1, \dots, \ell$. 
Then for regular $c$ the homogeneous component $\QH_r\c \subset \QH\c$ 
contains a copy of the defining module $V$ denoted by $V^{(r;c)}$ above in 
Section \ref{SSec:ReprQH}. Thus, the degree $r$ component $\QH_r$ of the 
space of universal quasiharmonics $\QH$ contains a unique copy $V^{(r)}$ of 
$\Complex[c]\otimes V$ (such that $\ev_c V^{(r)}= V^{(r;c)}$ for all 
regular $c$). On the other hand, the degree $r$ component 
$M_r^{(r/h)}(\triv)$ of the module $M^{(r/h)}(\triv)$ is known to contain a 
subspace $V_r$ of singular vectors isomorphic to $V$ (for the proof, see 
\cite{DunklJeuOpdam}, also \cite{BEG} and \cite{CE} for the case $W = A_n$, 
\cite{BEG} for $W = B_n$ and $W = D_n$, and \cite{Chmutova} for $W = 
I_2(m)$). We conjecture that this is a particular case of the following 
phenomenon:

 \begin{conjecture}\label{Cj:rhFlat}
Let $c = const. = r/h$ where $r \equiv d_i - 1 \bmod h$ for some $i = 1, 
\dots, \ell$. Then $\QH_r^{(r/h)} = \ev_{r/h}\QH_r$. In particular, $V_r = 
\ev_{r/h} V^{(r)}$.
 \end{conjecture}

We will prove  Conjecture \ref{Cj:rhFlat} for $W = I_2(m)$ in Section 
\ref{SSec:QH}. For the case $W = S_n$ see Section \ref{SSec:Comput}.

 \begin{definition} \label{Df:QuasiH}
We say that a positive integer $r$ is {\em $W$-good} if for each $c$ such 
that $\lmod c\rmod = r$ the Cherednik algebra $H_c(W)$ admits a 
finite-dimensional module of the form $S(V)/S(V) V_r$ where $V_r \subset 
S^r(V)$ is an irreducible $W$-submodule of dimension $\ell = \dim V$ 
(consequently, $V_r$ is a space of singular vectors, i.e.\ is killed by all 
the Dunkl operators).
 \end{definition}

 \begin{example}
According to \cite[Theorem 1.2]{BEG} (see also \cite{CE} for a shorter
proof) the number $r$ is $S_n$-good if $\name{gcd}(r,n) = 1$. By
\cite[Theorem 1.4]{BEG}, the number $r$ is $B_n$-good if it is odd and
$\name{gcd}(r,n) = 1$, and the number $r$ is $D_n$-good if it is odd and
$\name{gcd}(r,n-1) = 1$. According to \cite{Chmutova}, a number $r$ is
$I_2(m)$-good if it is not divisible by $m$ for $m$ odd, and by $m/2$ for
$m$ even.
 \end{example}

 \begin{conjecture}\label{Cj:123a}
Let $r$ be $W$-good. Then the quotient $S(V)/S(V)V^{(r;c)}$ is
finite-dimensional.
 \end{conjecture}

We will prove this conjecture for dihedral groups in Section 
\ref{SSec:Mult}.

 \begin{remark}
As we noted above (see beginning of Section \ref{SSec:Basic}), the ideal 
$S(V)V_r \subset S(V)$ is also a $H_{r/h}(W)$-submodule, so that the 
quotient $S(V)/S(V)V_r$ is a finite-dimensional $H_{r/h}(W)$-module. It is 
not true, though, for $S(V)/S(V)V^{(r;c)}$ in general, because the algebra 
$H_c(W)$ for regular $c$ has no finite-dimensional modules at all. 
Therefore each module over the $\Complex[c]$-algebra $H(W)$ is of infinite 
length over $\Complex[c]$ (i.e., after extending coefficients to the field 
$\Complex(c)$ it becomes an infinite-dimensional $\Complex(c)$-vector 
space).
 \end{remark}

\subsection{Quasiharmonics and elementary invariants for $W=S_n$}
\label{SSec:Comput}

In this section we make various observations concerning conjectures
\ref{Cj:SymmE}, \ref{Cj:rhFlat} and \ref{Cj:123a} in the case when $W = S_n$ is the symmetric group.

The reflections in $S_n$ are transpositions $(ij)$; they form a single conjugacy class. Therefore a conjugation-invariant function $c$ 
must be a constant. The defining module for $W = S_n$ is $V = 
\Complex^{n-1}$. It is convenient to assume that $V = \Complex^n/\Complex$ 
where the additive group $\Complex$ acts on $\Complex^n$ by simultaneous 
translations: $(x_1, \dots, x_n) \mapsto (x_1+b, \dots, x_n+b)$. The group 
$S_n$ acts in $V$ by permutation of variables. The symmetric algebra $S(V)$ 
is isomorphic to the algebra of all translation-invariant polynomials in 
$n$ variables: $f(x_1+b, \dots, x_n+b) = f(x_1, \dots, x_n)$ for all $b \in 
\Complex$.

Consider now the space $V^{(r;c)} \subset \QH_r\c(S_n)$ (where $n$ does not 
divide $r$), and let $q_{n,r}\c \in V^{(r;c)}$ be a non-zero invariant of 
the standard subgroup $S_{n-1} \subset S_n$. Then denote 
$q_{i,r}\c\bydef w(q_{1,r}\c)$ for any $w\in W$ such that $w(n)=i$. Clearly, the 
polynomials $q_{i,r}\c$ are well-defined.  

 \begin{proposition} 
For each $r>0$  there exists a constant $\alpha_r(c)$ such that for 
$i = 1,2,\dots,n$ one has:
 \begin{equation}\label{Eq:q-recursion}
\nabla_{x_i}\c q_{i,r}\c = \alpha_r(c) q_{i,r-1}\c,
 \end{equation}
where we use the convention $q_{i,kn}\c \bydef e_{kn}\c$ (see Corollary 
\ref{cor:multiple invariants}).
 \end{proposition}

 \begin{proof}
Let $\phi_{i,r-1}\c \bydef  \nabla_{x_i}\c q_{i,r}\c$ for all $r > 0$, 
$i=1, \dots, n$. Clearly, $\phi_{i,r}\c \in \QH_{r-1}\c$ and 
$w(\phi_{i,r}\c) = \phi_{w(i),r}\c$ for any $w \in S_n$. The sum 
$\sum_{i=1}^n \phi_{i,r}\c$ is a $S_n$-invariant element of $\QH_{r-1}\c$; 
so, by Corollary \ref{cor:multiple invariants}, if $r-1$ is not divisible by $n$, then 
the sum is zero. Thus, for $r \not\equiv 0,1 \bmod n$ the polynomials 
$\phi_{i,r}\c$, $i = 1, \dots, n$, form a copy of the defining $S_n$-module 
inside $\QH_{r-1}\c$. This copy must be $V^{(r-1,c)}$ by \eqref{Eq:MutDef}, 
so that \eqref{Eq:q-recursion} holds for all $r \not\equiv 0,1 \bmod n$. 

Let now $r = kn$ for some $k > 0$. Then each $q_{i,kn}\c$ is the 
quasiharmonic invariant $e_{kn}\c$, and the span of all $\phi_{i,kn-1}\c$ 
is isomorphic to $V$. Therefore, it is $V^{(kn-1,c)}$. Hence
 \begin{equation*}
\nabla_{x_i}\c e_{kn}\c = \alpha_{kn}(c) q_{i,kn-1}\c
 \end{equation*}
which proves \eqref{Eq:q-recursion} in this case.

Take now $r = kn + 1$. The elements $p_{kn,i}\c \bydef \nabla_{x_i}\c
q_{i,kn+1}\c$, $i = 1,\dots,n$ span an $S_n$-submodule $U_{kn}$ in
$\QH_{kn}\c$. This submodule is a quotient of the permutation $S_n$-module
$\triv \oplus V$. By \eqref{Eq:MutDef}, the $S_n$-module $\QH_{kn}\c$
contains no irreducible summands isomorphic to $V$  and, on the other hand,
$p_{kn,i}\c \ne 0$ for $c$ generic. Therefore, the $S_n$-module $U_{kn}$
must be trivial and $p_{kn,i}\c=\CC\cdot q_{kn}\c$ for $i = 1, \dots, n$.
That is,
 \begin{equation*}
\nabla_{x_i}\c q_{i,kn+1}\c = \alpha_{kn+1}(c) e_{kn}\c
 \end{equation*}
for  all $i = 1, \dots, n$.

The proposition is proved.
 \end{proof}

 \begin{remark}
The exact value of $\alpha_r(c)$ in \eqref{Eq:q-recursion} depends on the 
normalization of the polynomials $q_{i,r}\c$ (recall that $q_{i,r}\c$ are 
defined up to a multiplicative constant); explicit calculations for $n, r 
\le 5$ show that one can choose the normalization so that $\alpha_r(c) = r 
- nc$ and $q_{i,r}^{(r/n)} \not\equiv 0$. This supports Conjecture 
\ref{Cj:SymmE}.
 \end{remark}

 \begin{remark}
In the case $W = S_n$  the singular polynomials of $H_c(W)$ are known. If 
$c = r/n$ then the space $V_r$ of singular vectors in $S^r(V)$ is spanned 
by $f_{i,r}^{(r/n)}(x)$, $i = 1, \dots, n$, where $f_{i,r}^{(c)}(x)$ is the 
non-symmetric Jack polynomials defined by the formula
 \begin{equation*}
f_{i,r}\c(x) = \name{Res}_{z = \infty} \bigl( (1-x_1/z) \dots (1-x_n/z)
\bigr)^c \frac{z^rdz}{z - x_i}
 \end{equation*} 
(see \cite{DuInter} for proof).  Singular vectors for other $c \in \Sing$ 
are also specializations of some other non-symmetric Jack polynomials (see 
\cite{DuSing}, \cite{DuSingMod} for details). 

Note, however, that $q_{i,r}\c$ is not a multiple of $f_{i,r}\c$ for 
regular $c$. For instance, for $c = 0$ one has $f_{i,r}^{(0)}=x_i^r$, 
which, clearly, is not quasiharmonic. However, if $c = r/n$, Conjecture 
\ref{Cj:rhFlat} implies that $q_{i,r}\c$ and $f_{i,r}\c$  are proportional.
 \end{remark}

Conjecture \ref{Cj:rhFlat} claims that the family of vector spaces 
$\QH_r\c$ is flat over $c = r/h$. We prove in Section \ref{SSec:QH} 
(Proposition \ref{Pr:Dim}) that for the dihedral groups this family is flat 
everywhere. For $W = S_n$ it is not the case. Consider for example the case 
$n = 4$ ($S_3$ is the dihedral group $I_2(3)$). Computations using 
{\ttfamily MuPad} computer algebra system show that the space $\QH_3$ for 
the algebra $H(S_4)$ is isomorphic, as an $S_4$-module, to $V \oplus (V 
\otimes \eps)$ where $V$ is the defining module and $\eps$ is the sign 
character --- so, $\dim \QH_3\c = 6$ for generic $c$, as predicted by 
Proposition \ref{Pr:deformed hilbert series}. Exceptional values are $c = 
1/2$ and $c = 1/4$: here $\dim \QH_3\c = 7$; as a $S_4$-module one has then 
$\QH_3\c = V \oplus (V \otimes \eps) \oplus \triv$. The isotypic component 
of $V \subset \QH_3\c$ is $V^{(3;c)}$; the Dunkl operators $\nabla_{x_i}$, 
$i = 1, \dots, 4$, map it to zero for $c = 3/4$ confirming Conjecture 
\ref{Cj:SymmE}.

For $r = 4$ similar computations show that $\dim \QH_4(S_4) = 6$. As a
$S_4$-module $\QH_4\c$ for generic $c$ is isomorphic to $\triv \oplus \tau
\oplus (V \otimes \eps)$ where $\tau$ is the $2$-dimensional irreducible
module. Exceptional values are $c = 1/3$, $c = 1/2$ and $c = 3/4$ where
$\dim \QH_4\c = 9$ and $\QH_4\c = \triv \oplus \tau \oplus (V \otimes \eps)
\oplus V$.

Here are some formulas for the polynomials $e\c_k$ in the case $W = S_n$.

Recall that for $W = S_n$ the module $M(\triv) = S(V)$ is generated by
elements $x_1, \dots, x_n$ with the relation $x_1 + \dots + x_n = 0$.
Explicit computations show that
 \begin{align*}
&e\c_2 = e_2, \quad e\c_3 = e_3,\\
&e\c_4 = (n-2)(n-3)(1-nc) e_2^2/2 + (n^2(n-1)c - n(n+1)) e_4,\\
&e\c_5 = (n-3)(n-4)(1-nc) e_2 e_3 + (n^2(n-1)c - n(n+5)) e_5
 \end{align*}
where $e_s$ is the $s$-th elementary symmetric function.

Invariant quasiharmonic polynomials $e_{kh}\c$ (from Corollary 
\ref{cor:multiple invariants}) are computed by now only for the case $W = 
S_4$ and small $k$ (plus for all dihedral groups, see Section 
\ref{SSec:QH}, page \pageref{Pg:Qkh}). For example, one has
 \begin{align*}
&e_4\c = 4(12c - 5)e_4 - (4c - 1)e_2^2, \\
&e_8\c = (16c^2 - 32c + 27)e_2^4 - 24(16c^2 - 40c + 29) e_2^2 e_4 - 24(12c
- 13)e_2 e_3^2 \\
&\hphantom{e_8\c = (16c^2 - 32c + 27)e_2^4}+ 48(12c-13)(4c-5) e_4^2.
 \end{align*}
The authors are planning to write a separate paper dealing with properties 
of polynomials $e\c_k$ and $e_{kh}\c$.

We finish the section with an easy observation that one cannot replace
$\QH\c$ in conjectures \ref{Cj:SymmE} and \ref{Cj:rhFlat} by the space of
harmonics $\HH\c$. Namely, take $r = h+1$. For any irreducible reflection
group $W$ the smallest exponent is $d_1 = 2$, so that $r \equiv d_1-1 \bmod
h$. Then for regular $c$ one has $[V:\QH_r\c] = 1$ by \eqref{Eq:MutDef} and
$[V:\HH_r\c] = 0$ by  \eqref{Eq:standard multiplicity}. Thus, the space
$V^{(h+1;c)} \subset \QH_{h+1}\c$ is defined and not annihilated, for
regular $c$, by the operator $\nabla_{e_h}$. For $c = (h+1)/h \in \Sing$
the component $S^{h+1}(V)$ contains a subspace $V_{h+1}$ of singular
vectors isomorphic to the defining module $V$. Naturally, $V_{h+1} \subset
\HH_{h+1}^{(h+1)/h}$. So $0 = \ev_{(h+1)/h} \HH_{h+1} \ne
\HH_{h+1}^{(h+1)/h}$, i.e.\ the space of singular vectors cannot be
deformed in the class of harmonics. Conjecture \ref{Cj:SymmE} claims that for
quasiharmonics such situation is impossible --- they are, so to say,
flexible enough.

\section{Quasiharmonic elements: the dihedral group case}\label{Sec:Dih}

In this section we study quasiharmonic polynomials in the rational
Cherednik algebras for the dihedral group $W = I_2(m)$. In particular, we
verify conjectures \ref{Cj:SymmE}, \ref{Cj:rhFlat} (Corollary
\ref{Cr:Cj123valid}, Corollary \ref{Cr:Cj123validEven}) and \ref{Cj:123a}
(Corollary \ref{Cr:Cj123avalid}) for them. We also study the quotient
algebra of $M(\triv)$ by the submodule generated by homogeneous components
of $\QH\c$ (Section \ref{SSec:Mult}).

\subsection{Structure of the dihedral group (a summary)}\label{SSec:DihStr}

As an abstract group, $I_2(m)$ is generated by two elements $s_0, s_1$ with
the relations $s_0^2 = s_1^2 = (s_0 s_1)^m = 1$. Also, $I_2(m)$ is a finite
reflection group acting in the space $V = \Complex^2$ equipped with the
non-degenerate bilinear form $(\cdot,\cdot)$. The set $\Refl$ contains $m$
reflections $s_0, \dots, s_{m-1}$. The action of $I_2(m)$ is a
complexification of the action in $\Real^2$; to mark this fact we will be
using the basis $z \bydef x_1 + ix_2$, $\bar z \bydef x_1 - ix_2$ in
$\Complex^2$ instead of the usual $x_1, x_2$. Then the reflection $s_j$
acts as follows:
 \begin{equation}\label{Eq:Action}
s_j(z) = -\zeta^j \bar z, \quad s_j(\bar z) = -\zeta^{-j} z.
 \end{equation}
where $\zeta = e^{2\pi i/m}$ is the $m$-th primitive root of unity.

If $m$ is odd, then all the reflections $s \in \Refl$ are conjugate to one
another; so, $c = const.$ in the definition of the Cherednik algebra. If
$m$ is even, then the reflections $s \in \Refl$ split into two conjugacy
classes: $s_j$ with $j$ even and $s_j$ with $j$ odd. We will denote $c(s_j)
= c_1$ for the even class, and $c(s_j) = c_2$ for the odd class.

The $H_c(I_2(m))$-module $M(\triv) = S(V) = \Complex[z,\bar z]$ will be the
main object of study throughout Section \ref{Sec:Dih}.

The irreducible complex modules over the group $I_2(m)$ have dimensions $1$
and $2$ and can be described as follows:
 \begin{enumerate}
\item If $m$ is odd, then the group $I_2(m)$ has two $1$-dimensional
modules: the trivial one $\triv$ and the sign
one $\eps$ where $s_0 = s_1 = -1$.

\item If $m$ is even, the group has four $1$-dimensional representations:
$\triv$, $\eps$, $\mu_1$ and $\mu_2 = \mu_1 \otimes \eps$ where in $\mu_1$
one has $s_0 = 1$ and $s_1 = -1$, and in $\mu_2$, vice versa.

\item The $2$-dimensional irreducible $I_2(m)$-modules $Z_k$, $1 \le k < m/2$. The element $s_0$ acts on $Z_k\cong \Complex^2$ as multiplication by  $\left(\begin{array}{cc}0 & -1 \\ -1 & 0\end{array}\right)$, and
the element $s_1$ -- by 
$\left(\begin{array}{cc}   0 & -\zeta^{-k}\\ -\zeta^k & 0 \end{array}\right)$.
 \end{enumerate}

One can define the modules $Z_k$ for all $k \in \Integer$ by the
same formulas. It is easy to see that $Z_k\cong Z_{-k} \cong Z_{m+k}$; also $Z_0 = \triv \oplus \eps$ and, for $m$ even, $Z_{m/2} =
\mu_1 \oplus \mu_2$. Any irreducible $Z_r$ is isomorphic to a $Z_k$ with $1
\le k < m/2$. In particular, $Z_1$ is the defining $I_2(m)$-module.

This description allows to relate the modules $M(\tau)$ over the Cherednik
algebra $H_c(W)$ for different $1$-dimensional $\tau$. Thus, Dunkl
operators for $M(\eps)$ are the same as Dunkl operators for $M(\triv)$,
with the change $c \mapsto -c$. If $m$ is even, then the Dunkl operators
for $M(\mu_1)$ are the same as for $M(\triv)$, with the change $(c_1, c_2)
\mapsto (c_1, -c_2)$. By this reason, in the following sections we will
describe Dunkl and other related operators mostly for $M(\triv)$.

\subsection{Summary of main results}

Here we list, for the reader's convenience, the main results about dihedral
group case to be proved later in this section.

The first group of results are {\em proofs of some conjectures} mentioned
in Section \ref{Sec:Deform}:

 \begin{itemize}
\item Conjecture \ref{Cj:SymmE} is proved for the dihedral group $W =
I_2(m)$ with any $m$ and $c = const.$ (Corollary \ref{Cr:SymmEvalid}) and
for $m$ even and a generic $2$-valued $c$ (Corollary
\ref{Cr:SymmEvalidEven}).

\item Conjecture \ref{Cj:rhFlat} is proved for all dihedral groups
(Corollary \ref{Cr:Cj123valid}).

\item Conjecture \ref{Cj:123a} is proved for all dihedral groups (Corollary
\ref{Cr:Cj123avalid}).
 \end{itemize}

The second group of results are various {\em explicit formulas}:

 \begin{itemize}
\item In Section \ref{SSec:DuExpl} we derive formulas for polynomials
$e_{d_i}\c$ mentioned in Theorem \ref{Th:BetaIndep}. The exponents of the
dihedral group $I_2(m)$ are $2$ and $m$; formulas for $e_2\c$ and $e_m\c$
are given in  Corollary \ref{Cr:Beta} for $c = const.$ and in Corollary
\ref{Cr:Beta2C} for $m$ even and $c \ne const$.

\item In Section \ref{SSec:QH} we give formulas for quasiharmonic
polynomials of the dihedral group. According to Definition \ref{Df:QuasiH}
these polynomials are elements of the kernel of a single operator
$\nabla_{e_2}$. The dimension of this kernel is $2$ for every degree
(Proposition \ref{Pr:Dim}). Explicitly the basic quasiharmonics are given by
equations \eqref{Eq:ExplR}, \eqref{Eq:ExplRho} (see also \eqref{Eq:RRes}
and \eqref{Eq:NextR}); for the proof, see Theorem \ref{Th:DefR}.

\item In Section \ref{SSec:Mult} we study the quotient of the module
$M(\triv)$ by the ideal generated by quasiharmonics of some fixed degree
$n$. This quotient is a standard Frobenius algebra (see Appendix for
necessary definitions), and thus can be described in terms of its
characteristic polynomial. In Theorem \ref{Th:ChPoly} this polynomial is
computed explicitly for $c = const$.
 \end{itemize}

\subsection{Explicit formulas for Dunkl operators}\label{SSec:DuExpl}

In this section we will describe the Dunkl operators and the invariant
operator $\nabla_{e_2} = \nabla_{z \bar z}$ of degree $2$ in the module
$M(\triv)$ for the dihedral group $I_2(m)$. According to the remark in the
end of the previous section, it will allow us to obtain similar formulas for
$M(\tau)$ with any $1$-dimensional $\tau$. Formulas for $M(Z_k)$ also exist
but are cumbersome and do not serve our primary purpose (the study of
quasiharmonic elements), so we omit them.

The importance of the operator $\nabla_{e_2}$ comes from two reasons.
First, the dihedral group $I_2(m)$ has only two exponents: $d_1 = 2$ and $h
= d_2 = m$, so the space of quasiharmonics is simply the kernel of
$\nabla_{e_2}$. Second, there holds the following proposition:

 \begin{proposition}[\cite{BEG0}] \label{Pr:SL2} Let $E,F,H$ be endomorphisms of  $\Complex[z,\bar z]$ 
 given respectively by: $E = z \bar z $, $F= -\nabla_{e_2}$ and $H(P) = z \frac{\partial}{\partial z} + \bar z \frac{\partial}{\partial \bar z} +
(1 - mc)$. Then the operators $E$, $F$, and $H$ form a representation of the
Lie algebra $\name{\mathfrak{sl}}_2$, i.e.\ satisfy the relations
 \begin{equation*}
[E,F]=H, \quad [H,E]=2E, \quad [H,F]=-2F.
 \end{equation*}
 \end{proposition}
\noindent Below we will denote $\nabla_{e_2}$ by $-F$.

We consider two separate cases: the case $c = const.$ (possible for any
$m$) and the general case (possible only for $m$ even) when $c$ takes two
different values, $c_1$ and $c_2$, on even-numbered and odd-numbered
reflections.

Denote by $Y$, $\bar Y$ the Dunkl operators corresponding, in the notation
of \eqref{Eq:DunklOrig}, to the vectors $y = (1/2,-i/2)$ and $y = (1/2,
i/2)$, respectively. One can write symbolically $Y = \nabla_z$ and $\bar Y
= \nabla_{\bar z}$ because for $c = 0$ one has $Y = \partial/\partial z$
and $\bar Y = \partial/\partial \bar z$. Also let $T: \Complex^2 \to
\Complex^2$ denote the operator exchanging $z$ and $\bar z$, so that
$(TP)(z,\bar z) \bydef P(\bar z,z)$. Clearly,
 \begin{equation} \label{Eq:Exch}
Y T = T \bar Y, \qquad FT = TF.
 \end{equation}

 \begin{proposition} \label{Pr:DunklExpl}
Let $c = const$. Then in the module $M(\triv)$ one has
 \begin{equation*}
Y(z^a \bar z^b) = \begin{cases}
a z^{a-1} \bar z^b - mc \sum\limits_{0 \le k \le (a-b-1)/m} (-1)^{mk}
z^{a-mk-1} \bar z^{b+mk} &\text{if $a \ge b \ge 0$,} \\
a z^{a-1} \bar z^b + mc \sum\limits_{1 \le k \le (b-a)/m} (-1)^{mk}
z^{a+mk-1} \bar z^{b-mk} &\text{if $b \ge a \ge 0$.}
 \end{cases}
 \end{equation*}
 \end{proposition}

 \begin{proof}
Both equations are proved in a similar manner, so we concentrate on the
first one. Its first summand is equal to $\pder{z^a \bar z^b}{z}$;
therefore it is enough to prove that the sum in the right-hand side is
equal to $-cD_z(z^a \bar z^b)$ where $D_z$ is the ``difference term'' (the
sum over $s \in \Refl$) of the Dunkl operator, see \eqref{Eq:Dunkl}.

The polynomial $z \bar z$ is $I_2(m)$-invariant (it is the
elementary invariant $e_2$). So, if $a \ge b$, then $D_z(z^a \bar z^b) = (z
\bar z)^b D_z(z^{a-b})$. It is thus enough to consider the case $b = 0$, $a
\ge 0$.

Taking into account that  $z - s_j(z) 
= z + \zeta^j \bar z = \zeta^{j/2} \alpha_j(x)$  by \eqref{Eq:Action}, we obtain
$$D_z(z^a) = \sum_{j=0}^{m-1} \frac{z^a - s_j(\zeta^a)} {z -
s_j (z)}=\sum_{j=0}^{m-1} \frac{z^a - (-\zeta^{j} \bar z)^a}{z -
(-\zeta^{j} \bar z)} $$
$$= \sum_{j=0}^{m-1} \sum_{\ell=0}^{a-1} (-\zeta^{j} \bar
z)^\ell z^{a-1-\ell}= \sum_{\ell=0}^{a-1} (-1)^\ell z^{a-1-\ell} \bar z^\ell
\sum_{j=0}^{m-1} \zeta^{j\ell},
$$
and the identity
 \begin{equation*}
\sum_{j=0}^{m-1} \zeta^{j\ell} = \begin{cases}
m, &\text{if $m \divides \ell$},\\
0, &\text{otherwise}
\end{cases}
 \end{equation*}
finishes the proof.
 \end{proof}

Explicit formulas for $\bar Y$ can be obtained from Proposition
\ref{Pr:DunklExpl} using \eqref{Eq:Exch}.

 \begin{proposition} \label{Pr:FExpl}
Let $c = const$. Then in the module $M(\triv)$ for all $a \ge b \ge 0$ one
has
 \begin{equation}\label{Eq:FExpl}
F(z^a\bar z^b) = (mc-a)b z^{a-1} \bar z^{b-1} - mc\sum_{1 \le k \le
(a-b)/m} (-1)^{mk} (a-b-mk) z^{a-mk-1} \bar z^{b+mk-1}.
 \end{equation}
 \end{proposition}

 \begin{proof}
In the notation of Proposition \ref{Pr:DunklExpl} $Y = \pder{}{z} - cD_z$,
$\bar Y = \pder{}{\bar z} - cD_{\bar z}$, and therefore
 \begin{equation} \label{Eq:FViaD}
F = -\pdertwo{}{z}{\bar z} + c \bigl(\pder{}{\bar z} D_z + D_{\bar z}
\pder{}{z}\bigr) - c^2 D_{\bar z} D_z.
 \end{equation}

For $a = b$ the result is evident. If $a > b$ then one has
$-\pdertwo{}{z}{\bar z}(z^a \bar z^b) = -ab z^{a-1} \bar z^{b-1}$ and also,
by Proposition \ref{Pr:DunklExpl},
 \begin{align*}
\pder{}{\bar z} D_z (z^a \bar z^b) &= m\pder{}{\bar z}
\sum_{0 \le k \le (a-b-1)/m} (-1)^{mk} z^{a-mk-1} \bar z^{b+mk} \\
&= \sum_{0 \le k \le (a-b-1)/m} (b+mk) (-1)^{mk} z^{a-mk-1} \bar z^{b+mk-1}
\\
&= mb z^{a-1} \bar z^{b-1} + \sum_{1 \le k \le (a-b-1)/m} (b+mk) (-1)^{mk}
z^{a-mk-1} \bar z^{b+mk-1}
 \end{align*}
and
 \begin{equation*}
D_{\bar z} \pder{}{z} (z^a \bar z^b) = a D_{\bar z} (z^{a-1} \bar z^b) =
-ma \sum_{1 \le k \le (a-b-1)/m} (-1)^{mk} z^{a-mk-1} \bar z^{b+mk-1}.
 \end{equation*}
Thus, application to $z^a \bar z^b$ of all the terms of \eqref{Eq:FViaD}
except the last one already gives the right-hand side of \eqref{Eq:FExpl}.
As for the last term, write $a = mu + p$ where $u$ and $p$ are integers and
$0 \le p < m$. Then Proposition \ref{Pr:DunklExpl} implies that
 \begin{equation*}
D_{\bar z} (z^{mu + p}) = -m \sum_{k=1}^u z^{m(u-k)+p} \bar z^{mk-1} =
-D_{\bar z} (z^p \bar z^{mu}).
 \end{equation*}
It was noted (see the proof of Proposition \ref{Pr:DunklExpl}) that
$D_{\bar z}$ (as well as $D_z$) commute with the multiplication by $z \bar
z$. Therefore
 \begin{align*}
D_{\bar z} D_z(z^a) &= D_{\bar z} \sum_{k=0}^{[u/2]} (-1)^{mk}
(z\bar z)^{mk} (z^{m(u-2k)+p} + z^p \bar z^{m(u-2k)}) \\
&= \sum_{k=0}^{[u/2]} (-1)^{mk} (z\bar z)^{mk} D_{\bar z} (z^{m(u-2k)+p} +
z^p \bar z^{m(u-2k)}) = 0,
 \end{align*}
and so $D_{\bar z} D_z (z^a \bar z^b) = (z \bar z)^b D_{\bar z} D_z
(z^{a-b}) = 0$ for all $a > b$. This completes the proof of
\eqref{Eq:FExpl}.
 \end{proof}
\noindent Explicit formulas for $F(z^a \bar z^b)$ with $a \le b$ can be
obtained using \eqref{Eq:Exch}.

 \begin{corollary} \label{Cr:Beta}
For $c = const.$ the polynomial $z^m + \bar z^m$ is quasiharmonic. In terms
of Proposition \ref{Pr:unique elementary invariant} one has $\eps_2 = z \bar z$ and
$\eps_m = z^m + \bar z^m$. These polynomials are algebraically independent
for all $c$.
 \end{corollary}

Consider now the case of general $c$ for $m$ even, so that $c(s_j) = c_1$
for $j$ even and $c(s_j) = c_2$ for $j$ odd.

 \begin{proposition} \label{Pr:DunklExpl2C}
In the module $M(\triv)$ with $m$ even for any $a \ge b \ge 0$ one has
 \begin{equation*}
Y(z^a \bar z^b) = az^{a-1} \bar z^b - \frac{m}{2}\sum_{0 \le k \le
2(a-b-1)/m} (-1)^{mk/2} z^{a-1-mk/2}\bar z^{b+mk/2} ((-1)^k c_1
+ c_2).
 \end{equation*}
For any $b \ge a \ge 0$ one has
 \begin{equation*}
Y(z^a \bar z^b) = az^{a-1} \bar z^b + \frac{m}{2}\sum_{1 \le k \le
2(b-a)/m} (-1)^{mk/2} z^{a-1+mk/2}\bar z^{b-mk/2} ((-1)^k c_1
+ c_2).
 \end{equation*}
 \end{proposition}

 \begin{proposition} \label{Pr:FExpl2C}
In the module $M(\triv)$ with $m$ even for any $a \ge b \ge 0$ one has
 \begin{multline*}
F(z^a \bar z^b) = (m(c_1 + c_2)/2 - a)b z^{a-1} \bar z^{b-1} \\
- \frac{m}{2} \sum_{1 \le k \le 2(a-b)/m} (-1)^{mk/2} z^{a-1-mk/2} \bar
z^{b-1+mk/2} (a-b-mk/2) ((-1)^k c_1 + c_2).
 \end{multline*}
 \end{proposition}

Proofs of these propositions are similar to those of Propositions
\ref{Pr:DunklExpl} and \ref{Pr:FExpl}. Explicit formulas for the operator
$\bar Y$  and for $F(z^a \bar z^b)$ with $b \ge a \ge 0$ can be obtained
using \eqref{Eq:Exch}.

 \begin{corollary} \label{Cr:Beta2C}
For $m$ even the polynomial
 \begin{equation*}
e_m\c = (c_1+c_2-1)(z^m + \bar z^m) + 2(-1)^{m/2}(c_2-c_1)(z \bar z)^{m/2}
 \end{equation*}
is quasiharmonic. In terms of Proposition \ref{Pr:unique elementary invariant}
(and Remark \ref{Rm:DefEC}) one has $e_2\c = z \bar z$. These polynomials
are algebraically independent for all the functions $c$ such that $c_1 +
c_2 \ne 1$.
 \end{corollary}

\subsection{Quasiharmonic polynomials}\label{SSec:QH}

In this section we describe the space $\QH_n\c \subset M(\triv)$, i.e.\ the
kernel of the operator $F$ described above.

 \begin{proposition} \label{Pr:Dim}
For all $c$ and $n$ the space $\QH_n\c$ has dimension $2$ and is
isomorphic, as an $I_2(m)$-module, to $Z_n$.
 \end{proposition}

 \begin{proof}
Consider the action of $\name{\mathfrak{sl}_2}$ on $M(\triv) = \Complex[z,
\bar z]$ described in Proposition \ref{Pr:SL2}. By a classical theorem,
$\Complex[z, \bar z]$ splits into a direct sum of spaces spanned by the
bases $a, Ea, E^2a, \dots$, where $\name{deg} a = r$ for some $r$ and
therefore $\name{deg} E^k a = r+2k$. On such a basis the operator $F$ acts
as $FE^k a = 2k E^{k-1} a$. The operator $E$ has a trivial kernel;
therefore the operator $F: \Complex_n[z, \bar z] \to \Complex_{n-2}[z, \bar
z]$ has a trivial co-kernel. This implies the equality $\dim \QH_n\c = \dim
\Complex_n[z, \bar z] - \dim \Complex_{n-2}[z, \bar z] = 2$ for all $c$ and
$n$.

The space $M_n(\triv) = \Complex_n[z,\bar z]$ splits into a sum of
$I_2(m)$-isotypic components corresponding to irreducible representations
described in Section \ref{SSec:DihStr}.
Note now that the space $\QH_n\c \subset M_n(\triv)$ for any $n$ and $c$
lies in the kernel of the operator $F$ and therefore cannot lie in the
image of $E$.

Thus, $\QH_n\c$ contains polynomials $P_1 = z^n + \dots$ and $P_2 = \bar
z^n + \dots$ where dots mean a linear combination of $z^{n-k} \bar z^k$
with $k = 1, \dots, n-1$. Since $\dim \QH_n\c = 2$, it is isomorphic to the
$I_2(m)$-module generated by $z^n$ and $\bar z^n$, that is, to $Z_n$.
 \end{proof}

In particular, if $n \equiv \pm 1 \bmod m$, then $\QH_n\c$ is isomorphic to
the defining $W$-module.

Suppose now that the function $c$ is constant. Let $n = mq + r$ where $q,r$
are nonnegative integers and $0 \le r \le m-1$. For a nonnegative integer
$k$ denote $\lambda_k(c) = c(c-1)\dots(c-k)$ (with the convention
$\lambda_{-1}(c) = 1$). Consider a polynomial
 \begin{equation}\label{Eq:ExplR}
R_{n,c}(z,\bar z) = \sum_{p = 0}^q (-1)^{mp} \binom{q}{p} \lambda_{q-p}(c)
\lambda_{p-1}(c) z^{n-mp} \bar z^{mp},
 \end{equation}
and take
 \begin{equation}\label{Eq:ExplRho}
\rho_{n,c}(z,\bar z) \bydef R_{n,c}/\lambda_{[q/2]}(c),
 \end{equation}

Equation \eqref{Eq:ExplR} can be rewritten in another form:
 \begin{equation} \label{Eq:RRes}
R_{n,c}(z,\bar z) = cq! z^r \name{Res\,}\limits_{w=0} (1 + (-\bar z)^m/w)^c
(1 + (-z)^m/w)^{c-1} w^{q-1}dw.
 \end{equation}
For $c = n/m$ this becomes a formula for singular vectors from
\cite{DunklJeuOpdam}.

One more form of \eqref{Eq:ExplR} is recursive:
 \begin{equation}\label{Eq:NextR}
R_{n,c} = \begin{cases}
z R_{n-1,c}, &r \ne 0,\\
(c-q) z R_{n-1,c} + (-1)^{mq} c \bar z \bar R_{n-1,c}, &r = 0.
 \end{cases}
 \end{equation}
Proofs of \eqref{Eq:RRes} and \eqref{Eq:NextR} are immediate.

 \begin{proposition}[\cite{DunklDef}, p.~181] \label{Pr:RQuasiH}
Suppose that $c = const$. Then the polynomial $\rho_{n,c}$ is quasiharmonic
for all $n$ and $c$. One has
 \begin{equation*}
Y(\rho_{n,c}) = \begin{cases}
(n-mc) \rho_{n-1,c} &\text{if\ } r \ne 0,\\
mq(2c-q) \rho_{n-1,c} &\text{if $r = 0$ and $q$ is odd},\\
2mq \rho_{n-1,c} &\text{if $r = 0$ and $q$ is even},
\end{cases}
 \end{equation*}
and
 \begin{equation*}
\bar Y (\rho_{n,c}) = 0
 \end{equation*}
for all $n$ and $c$.
 \end{proposition}

Proposition \ref{Pr:RQuasiH} implies the following description of the space
of quasiharmonics with $c = const.$:

 \begin{theorem} \label{Th:DefR}
Suppose that $c = const$. If $r \ne 0$ then the polynomials $\rho_{n,c}$
and $\bar \rho_{n,c} \bydef T \rho_{n,c}$ form a basis for the space
$\QH_n\c$. If $r = 0$ (that is, $n = mq$) then the basis consists of two
polynomials $\rho_{mq,c} + \bar \rho_{mq,c}$ and $(\rho_{mq,c} - \bar
\rho_{mq,c})/(2c-q)$.
 \end{theorem}

 \begin{remark}
In fact, $\rho_{mq,c}$ and $\bar\rho_{mq,c}$ form a basis in $\QH_{mq}\c$,
too, unless $q$ is odd and $c = q/2$; see the proof below.
 \end{remark}

 \begin{proof}
In view of Proposition \ref{Pr:Dim} and Proposition \ref{Pr:RQuasiH} it
suffices to prove that $\rho_{n,c}$ and $\bar \rho_{n,c}$ are linearly
independent unless $r = 0$, $q$ is odd, and $c = q/2$. The element $\gamma
= s_0 s_1 \in W \subset H_c(W)$ (where $s_0, s_1$ are the defining
reflections) acts on $M(\triv) = \Complex[z,\bar z]$ as follows: $\gamma
P(z,\bar z) = P(\zeta z, \bar\zeta \bar z)$ where $\zeta = \exp(2\pi i/m)$
is a primitive $m$-th root of unity. The polynomials $\rho_{n,c}$ and $\bar
\rho_{n,c}$ are eigenvectors of the operator $\gamma$ with eigenvalues
$\zeta^n$ and $\zeta^{-n}$, respectively. If $r \ne 0$ then the eigenvalues
are different, and $\rho_{n,c}$ and $\bar \rho_{n,c}$ can be linearly
dependent only if they are identically zero, which never happens.

Let now $r = 0$, that is, $n = mq$. By \eqref{Eq:ExplR} and
\eqref{Eq:ExplRho} the coefficients at $z^{mk}\bar z^{m(q-k)}$ and
$z^{ml}\bar z^{m(q-l)}$ of the polynomials $\rho_{n,c}$ and $\bar
\rho_{n,c}$ form a matrix
 \begin{equation*}
A_{kl} = \left(\begin{array}{cc}
\frac{\lambda_{k-1}(c) \lambda_{q-k}(c)}{\lambda_{[q/2]}(c)} &
\frac{\lambda_{l-1}(c)\lambda_{q-l}(c)}{\lambda_{[q/2]}(c)} \\
\frac{\lambda_{q-k-1}(c) \lambda_k(c)}{\lambda_{[q/2]}(c)} &
\frac{\lambda_{q-l-1}(c) \lambda_l(c)}{\lambda_{[q/2]}(c)}
\end{array}\right).
 \end{equation*}
Apparently, $\rho_{n,c}$ and $\bar\rho_{n,c}$ are linearly dependent only
if $\det A_{k,l}(c) = 0$ for all $k$ and $l$.

One has
 \begin{align}
\det A_{kl} &= \frac{\lambda_{k-1}(c) \lambda_{l-1}(c)
\lambda_{q-k-1}(c) \lambda_{q-l-1}(c)}{\lambda_{[q/2]}^2(c)}
\det \left(\begin{array}{cc}
c-q+k & c-q+l \\
c-k & c-l
\end{array}\right) \nonumber\\
&= (k-l)(2c-q) \frac{\lambda_{k-1}(c) \lambda_{l-1}(c) \lambda_{q-k-1}(c)
\lambda_{q-l-1}(c)} {\lambda_{[q/2]}^2(c)}.\label{Eq:Akl}
 \end{align}

Suppose first that $q$ is even, so that $[q/2] = q/2$. Then $\det A_{kl}(c)$
can be zero only if $c$ is an integer between $0$ and $q$.
Now one has $\det A_{0q} = q(2c-q)
\frac{c\lambda_{q-1}(c)}{\lambda_{q/2}(c)^2}$, so that $\det A_{0q}(0) \ne
0$. Also $\det A_{01} = (2c-q) \frac{c \lambda_{q-1}(c)
\lambda_{q-2}(c)}{\lambda_{q/2}(c)^2}$, and therefore $\det A_{01}(c) \ne
0$ for $c = 1, \dots, q/2-1$. Similarly, $\det A_{q/2,q/2+1}(c) \ne 0$ for
$c = q/2, \dots, q$.

Thus, for every $c$ there exist $k,l$ such that $\det A_{kl}(c) \ne 0$, and
the polynomials $\rho_{n,c}$ and $\bar \rho_{n,c}$ are linearly
independent.

If $q$ is odd then the reasoning is the same except for the case $c = q/2$.
One gets $\det A_{kl}(q/2) = 0$ for all $k,l$, so that $\rho_{n,c}$ and
$\bar \rho_{n,c}$ for $c = q/2$ are linearly dependent (actually, equal). A
similar computation shows that in this case $\rho_{mq,c} + \bar
\rho_{mq,c}$ and $(\rho_{mq,c} - \bar \rho_{mq,c})/(2c-q)$ form a basis in
$\QH_{mq}\c$.
 \end{proof}

 \begin{corollary}[\cite{DunklJeuOpdam}] \label{Cr:DJO}
A function $c = const.$ is singular for the module $M(\triv)$ of the
dihedral group $W = I_2(m)$ if $c = n/m$, where $n$ is a positive integer
not divisible by $m$, or $c = \ell + 1/2$, where $\ell$ is a nonnegative
integer.
 \end{corollary}

 \begin{corollary}[of Proposition \ref{Pr:RQuasiH}, Theorem \ref{Th:DefR}
and Proposition \ref{Pr:Dim}]\label{Cr:Cj123valid}
Conjectures \ref{Cj:SymmE} and \ref{Cj:rhFlat} are valid for the dihedral
group $W = I_2(m)$ and $c = const$.
 \end{corollary}

 \begin{proof}
For $c = n/m$ the deformation in question is $V_r(c) = \QH_r\c$. For $c
= \ell + 1/2$ the deformation is $\rho_{m(2\ell+1),c}$.
 \end{proof}

\label{Pg:Qkh} Note that the polynomials $e_{km}\c = \rho_{km,c} + \bar \rho_{km,c}$, 
$k=1,2,\dots$ are exactly the quasiharmonic invariants mentioned in 
Corollary \ref{cor:multiple invariants}.  So, Theorem \ref{Th:DefR} and 
Proposition \ref{Pr:RQuasiH} imply

 \begin{corollary}[of Theorem \ref{Th:DefR} and Corollary
\ref{Cr:DJO}]\label{Cr:SymmEvalid}
Conjecture \ref{Cj:SymmE} is valid for dihedral groups with $c = const$.
 \end{corollary}

Later we are going to prove (see Corollary \ref{Cr:SymmEvalidEven}) a
similar statement for $I_2(m)$ with $m$ even and generic function $c$
assuming $2$ values.

Recall that  the space $\QH_n\c$ is isomorphic to $Z_n$ by Proposition
\ref{Pr:Dim}. Denote $S_{n,c} \in \QH_n\c$ an eigenvector with the
eigenvalue $\zeta^n$ of the operator $\gamma = s_0 s_1$. If $n$ is not divisible by $m/2$, then the
two eigenvalues of $\gamma$ are different and $S_{n,c}$ is defined uniquely
up to proportionality; we normalize it by the condition $S_{n,c} =
\lambda_{[2(n-1)/m]}(c_1+c_2) z^n + \dots$ (recall that we denote
$\lambda_k(c) = c(c-1)\dots(c-k)$; square brackets $[\alpha]$ mean the
biggest integer not exceeding $\alpha$). If $2n/m = q \in \Integer$, then
$\gamma = (-1)^q\name{Id}$; in this case we define $S_{n,c} = \sum_{k =
0}^{2q} r_k z^{n - mk/2} \bar z^{mk/2} \in \QH_n\c$ by the conditions $r_0
= \lambda_{[2(n-1)/m]}(c_1+c_2)$, $r_{2q} = 0$. Note that if $2n/m \notin
\Integer$, then it follows from Proposition \ref{Pr:FExpl2C} that $S_{n,c} =
\sum_{0 \le k \le 2n/m} r_k z^{n - mk/2} \bar z^{mk/2}$ does not contain
the term $\bar z^n$ either. So the coefficient at $\bar z^n$ in $S_{n,c}$
vanishes for all $n$, and the coefficient at $z^n$ is always equal to
$\lambda_{[2(n-1)/m]}(c_1+c_2)$. Also, denote $\bar S_{n,c} \bydef
TS_{n,c}$; then $\bar S_{n,c} = \lambda_{[2(n-1)/m]}(c_1+c_2)\bar z^n +
\dots$ and does not contain the term $z^n$.

Clearly, for $c$ generic the polynomials $S_{n,c}$ and $\bar S_{n,c}$
are linearly independent and therefore form a basis in $\QH_n\c$.

Polynomials $S_{n,c}$ (with a different normalization) were first
considered in \cite{DunklJeuOpdam}. The action of Dunkl operators on
$\QH_n\c$ is described as follows:

 \begin{proposition}[\cite{DunklJeuOpdam}] \label{Pr:ActDunklR2C}
Let $m$ be even, $n = mq/2 + r$, $q,r \in \Integer$, $0 \le r \le m/2-1$.
Then for generic $2$-valued function $c$ one has
 \begin{equation}\label{Eq:ActZ2C}
Y(S_{n,c}) = \begin{cases}
\left(n-\frac{m}{2}(c_1+c_2)\right) S_{n-1,c},
& r \ne 1,\\
 \begin{aligned}
\biggl(n-&\frac{m}{2}(c_1+c_2)\biggr) \biggl((c_1+c_2-q) S_{n-1,c}\\
&+ (-1)^{mq/2}((-1)^q c_1 + c_2)\bar S_{n-1,c} \biggr)
 \end{aligned}
&r = 1
 \end{cases}
 \end{equation}
and
 \begin{equation}\label{Eq:ActZBarR2C}
\bar Y(S_{n,c}) = \begin{cases}
0, & r \ne 0,\\
(-1)^{mq/2}\frac{m}{2}((-1)^q c_1 + c_2) \bar S_{n-1,c}, & r = 0.
 \end{cases}
 \end{equation}
 \end{proposition}

 \begin{corollary} \label{Cr:Cj123validEven}
Conjecture \ref{Cj:SymmE} is valid for dihedral groups $I_2(m)$ with $m$
even, and generic $c$.
 \end{corollary}

The polynomial $e_{mk}\c = S_{mk,c} + \bar S_{mk,c}$ is the quasiharmonic
invariant from Corollary \ref{cor:multiple invariants}.

 \begin{corollary}\label{Cr:SymmEvalidEven}
Conjecture \ref{Cj:SymmE} is valid for dihedral groups with even $m$ and
generic $2$-valued function $c$.
 \end{corollary}

 \begin{remark}
In view of Corollary \ref{Cr:SymmEvalidEven}, Conjecture \ref{Cj:SymmE} for
dihedral groups is reduced to the following statement: the kernel of Dunkl
operators applied to quasiharmonic invariants is bigger than usual exactly
when the Harish-Chandra pairing has a nontrivial kernel. We proved this for
$c = const.$ (see Corollary \ref{Cr:SymmEvalid} above).
 \end{remark}

Explicit (though rather complicated) formulas for the polynomials $S_{n,c}$ 
can be found in \cite{DunklJeuOpdam}. In particular, Proposition \ref{Pr:ActDunklR2C} 
implies the following recursive formula for $S_{n,c}$ similar to 
\eqref{Eq:NextR}:

 \begin{equation}
 \label{Eq:Recur2C}
S_{n+1,c} = \begin{cases}
z S_{n,c}, & r \ne 0,\\
(c_1+c_2-q) z S_{n,c} + (-1)^{mq/2} ((-1)^q c_1 + c_2) z \bar S_{n,c},
&r = 0.
\end{cases}
\end{equation}

\subsection{Multiplication modulo quasiharmonics}\label{SSec:Mult}

The module $M(\triv)$ is isomorphic to $\Complex[z,\bar z]$; so, it possess
a natural algebra structure. Denote by $J_{n,c} \subset \Complex[z,\bar z]$
the ideal generated by the space $\QH_n\c$ of quasiharmonics. We are
going to investigate the multiplication in the algebra ${\mathcal P}_{n,c}
= M(\triv)/J_{n+1,c}$.

We will be using some known facts from the theory of standard Frobenius
algebras(the algebra ${\mathcal P}_{n,c}$ is standard Frobenius for $c$ generic).

The necessary statements, as well as notation, are summarized in
the Appendix (Section \ref{Sec:Apolar}).

 \begin{lemma}\label{Lm:MultZZBar}
For every $c$ and for every $P \in J_{n,c}$ one has $z \bar z P \in
J_{n+1,c}$. If $c = const.$ is regular and $z \bar z P \in J_{n+1,c}$ then
$P \in J_{n,c}$.
 \end{lemma}

 \begin{proof}
The first statement follows obviously from  \eqref{Eq:NextR} with
$c = const.$ and from \eqref{Eq:Recur2C} for $m$ even and general
$c$.

Apply now the map $T$ (exchanging $z$ and $\bar z$) to \eqref{Eq:NextR}:
 \begin{equation*}
\bar R_{n,c} = \begin{cases}
\bar z \bar R_{n-1,c}, &r \ne 0,\\
(c-q) \bar z \bar R_{n-1,c} + (-1)^{mq} c z R_{n-1,c}, &r = 0.
 \end{cases}
 \end{equation*}
Using this together with \eqref{Eq:NextR} one obtains the relation
 \begin{equation*}
z\bar z R_{n-1,c} = \begin{cases}
\bar z R_{n,c}, &r \ne 0,\\
\frac{1}{q(q-2c)} (\bar z R_{n,c} - (-1)^n c z \bar R_{n,c}), &r = 0.
 \end{cases}
 \end{equation*}
which implies the second statement of the lemma.
 \end{proof}

 \begin{remark}
It is natural to expect that for $m$ even and regular $c \ne const.$ the
converse of Lemma \ref{Lm:MultZZBar} is also true. Numerical examples support this expectation
but for the moment of writing the proof was not known.
 \end{remark}

 \begin{proposition} \label{Pr:FinDim}
The algebra ${\mathcal P}_{n,c}$ is finite-dimensional for regular $c =
const.$ and for generic $c \ne const.$ if $m$ is even.
 \end{proposition}

 \begin{corollary}\label{Cr:Cj123avalid}
Conjecture \ref{Cj:123a} holds for the dihedral group $I_2(m)$.
 \end{corollary}

 \begin{proof}[of Proposition \ref{Pr:FinDim}]
We are going to prove that $M_{2n+1}(\triv) \subset J_{n+1,c}$ using
induction by $n$. The base $n = 1$ can be checked immediately. Suppose that
for $n-1$ the assertion is proved. It follows then from Lemma
\ref{Lm:MultZZBar} that $z\bar z M_{2n-1}(\triv) \subset J_{n+1,c}$. For $c
= 0$ the ideal $J_{n+1,0}$ is generated by the polynomials $z^{n+1}$ and
$\bar z^{n+1}$; therefore the generators $\rho_{n+1,c}$ (for $c = const.$)
and $S_{n+1,c}$ (for generic $c$ and $m$ even) have nonzero coefficient
$b(c)$ at $z^{n+1}$ for generic $c$.

If $c = const.$, then Theorem \ref{Th:DefR} and Corollary \ref{Cr:DJO}
imply a stronger statement: the coefficient $b(c)$ is nonzero for every
regular $c$.

Let now $n+1 = mq + r$. If $r \ne 0$ then for regular $c$ one has
$z^n\rho_{n+1,c} = b(c) z^{2n+1} + z \bar z \phi(c)$ with $b(c) \ne 0$, and
therefore $z^{2n+1} \in J_{n+1,c}$.

Similarly, $\bar z^{2n+1} \in J_{n+1,c}$, so that $M_{2n+1}(\triv) \in
J_{n+1,c}$ in this case.

If $r = 0$ (that is, $n+1 = mq$), consider the polynomial $\det A_{0q}(c)$
from the proof of Theorem \ref{Th:DefR}. It is nonzero for $c = 0$ and
therefore nonzero for generic $c$ ($c \ne const.$ included). Explicit
formulas \eqref{Eq:Akl} imply also that $A_{0q}(c) \ne 0$ for any regular
$c = const$. So, in both cases there exists a linear combination $\mu
\bydef \alpha z^n \rho_{n,c} + \beta \bar z^n \bar \rho_{n,c}$ such that
$\mu = b(c) z^{2n+1} + z\bar z \phi$ with $b(c) \ne 0$. So, $z^{2n+1} \in
J_{n+1,c}$ for $r = 0$, too. Similarly, $\bar z^{2n+1} \in J_{n+1,c}$, and
the induction is finished.
 \end{proof}

It follows now from Proposition \ref{Pr:AllMinim2} that the algebra
${\mathcal P}_{n,c}$ is standard Frobenius. Its degree is $2n$; denote
$p_{2n,c}$ its characteristic polynomial. It is defined up to a constant
multiple; we do not fix the normalization unless otherwise stated.

 \begin{proposition} \label{Pr:Laplace}
Let $c = const$. For every $n$ there exists a rational function $a_n(c)$
such that $\Delta p_{2n,c} = a_n(c) p_{2(n-1),c}$ (here $\Delta =
-\pdertwo{}{z}{\bar z}$ is the Laplace operator).
 \end{proposition}

 \begin{proof}
Let $u(z,\bar z)$ be a polynomial of degree $2n-2$. By the definition of
the characteristic polynomial, $u \in J_{n,c}$ if and only if $\langle
p_{2(n-1),c}, u\rangle = 0$. By Lemma \ref{Lm:MultZZBar} for generic $c =
const.$ this is equivalent to $z \bar z u \in J_{n+1,c}$, and therefore $0
= \langle p_{2n,c}, z \bar z u\rangle = \langle \Delta p_{2n, c},
u\rangle$. So, the annihilators of the polynomials $p_{2(n-1),c}$ and
$\Delta p_{2n,c}$ are the same, and therefore these polynomials are
proportional: $\Delta p_{2n,c} = a_n(c) p_{2(n-1),c}$. Here $a_n$ is a
rational function of $c$; if we choose the normalization of $p_{2n,c}$ for
all $n$ as in \eqref{Eq:QuasiRes}, then $a_n(c)$ will be actually a
polynomial.
 \end{proof}

For $c = const.$ we will compute $p_{2n,c}$. Denote $Q_1 = (c-1)z^{2m} +
(-1)^m c z^{m} \bar z^m$, and $Q_k = \prod_{i=0}^k (c-i) \cdot z^{2km} +
(-1)^{mk+m+k} \prod_{i=-1}^{k-1} (c+i) \cdot z^{m(k+1)} \bar z^{m(k-1)}$
for $k \ge 2$.

 \begin{lemma}
Let $n = mq + r$. Then in the notation of \eqref{Eq:ExplR} one has:
 \begin{align}
z^{mq} \bar z^r R_{n,c} &= \sum_{k=0}^{q-2} (-1)^{mk} \binom{q}{k}
\prod_{i=0}^{k-1} (c-i) \cdot (z \bar z)^{mk+r} Q_{q-k} \label{Eq:MultZ}\\
&\hphantom{=\,\,}
+ (-1)^{mq} \prod_{i=0}^{q-1} (c-i) \cdot (z \bar z)^{m(q-1)+r}
Q_1.\nonumber
 \end{align}
(an empty product for $k = 0$ is assumed to be $1$).

 \end{lemma}

 \begin{proof}
By \eqref{Eq:ExplR} it suffices to prove \eqref{Eq:MultZ} for $r = 0$. We
will be using the notation of \eqref{Eq:ExplR} where $\lambda_k(c) \bydef
c(c-1) \dots (c-k)$.

The term containing $Q_{r}$ in the right-hand side is actually equal to
 \begin{align*}
(-1)^{mk} &\lambda_{k-1}(c) (z \bar z)^{mk} (\lambda_{q-k}(c)
z^{2m(q-k)} \bar z^{mk} \\
&\hphantom{\lambda_{k-1}(c)}+ (-1)^{m(q-k)+m+q-k} \lambda_{q-k}(c+q-k-1)
z^{m(q-k+1)} \bar z^{m(q-k-1)}) \\
&= (-1)^{mk} \lambda_{k-1}(c) \lambda_{q-k}(c) z^{m(2q-k)} \bar
z^{mk}\\
&\hphantom{\lambda_{k-1}(c)}+ (-1)^{m(q+1)+q-k} \lambda_{k-1}(c)
\lambda_{q-k}(c+q-k-1) z^{m(q+1)} \bar z^{m(q-1)},
 \end{align*}
(where we abbreviated $r=q-k$) and the $Q_1$ term is
 \begin{align*}
(-1)^{m(q+1)} \lambda_0(c) \lambda_{q-1}(c) z^{mq} \bar z^{mq} + (-1)^{mq}
\lambda_0(c-1) \lambda_{q-1}(c) z^{m(q+1)} \bar z^{m(q-1)}.
 \end{align*}
So, every monomial $z^{m(2q-p)} \bar z^{mp}$ with $0 \le p \le q,\,p \ne 
q-1$ enters the right-hand side of \eqref{Eq:MultZ} exactly once, and the 
coefficient is the same as in the left-hand side (cf.\ \eqref{Eq:ExplR}).

The coefficient at $z^{m(q+1)} \bar z^{m(q-1)}$ in the left-hand side is
equal, by \eqref{Eq:ExplR}, to $(-1)^{m(q-1)} \binom{q}{q-1} \lambda_1(c)
\lambda_{q-2}(c) = (-1)^{m(q+1)} q c(c-1) \lambda_{q-2}(c)$, while the
corresponding coefficient in the right-hand side is
 \begin{align*}
\sum_{k=0}^{q-2} &(-1)^{m(q+1)+q-k} \binom{q}{k} \lambda_{k-1}(c)
\lambda_{q-k}(c+q-k-1) + (-1)^{m(q+1)} (c-1)\lambda_{q-1}(c) \\
&= (-1)^{m(q+1)}c(c-1) \sum_{0 \le k \le q, k \ne q-1} (-1)^{q-k}
\binom{q}{k} \lambda_{q-2}(c+q-k-1).
 \end{align*}
$\lambda_{q-2}(x)$ is a polynomial of degree $q-1$, so its $q$-th iterated
difference is identically zero: 
 \begin{equation*}
\sum_{k=0}^q (-1)^k \binom{q}{k} \lambda_{q-2}(x-k+1) = 0.
 \end{equation*}
Taking $x = c+q$, we obtain the desired identity for the coefficients of 
$z^{m(q+1)} \bar z^{m(q-1)}$. 
 \end{proof}

 \begin{proposition} \label{Pr:QinJ}
For any $c$ one has $z^r \bar z^r Q_q \in J_{n+1,c}$.
 \end{proposition}

 \begin{proof}
Proposition \ref{Pr:DerivDep} implies that the corollary is true when $c$
belongs to some closed subset of $\Complex$; therefore it is enough to
prove it for generic $c$.

Use induction in $q$. For $q = 1$ the inclusion can be checked immediately.
For bigger $q$, it follows from the induction hypothesis and Lemma
\ref{Lm:MultZZBar} that all the terms in \eqref{Eq:MultZ}, except $z^r \bar
z^r Q_q$, are in $J_{n,c}$ --- therefore, $z^r \bar z^r Q_q \in J_{n,c}$,
too.
 \end{proof}

 \begin{corollary} \label{Cr:Zeros}
If $c = -1, -2, \dots, -q+1$, then $z^{2n-r} \bar z^r \in J_{n+1,c}$.
 \end{corollary}

 \begin{proposition} \label{Pr:Monom}
If $c = 1, \dots, q$, then $z^{n+m+r} \bar z^{n-m-r} \in J_{n+1,c}$.
 \end{proposition}

 \begin{proof}
Similar to Proposition \ref{Pr:QinJ} one proves that $z^{2r}Q_q \in
J_{n+1,c}$ for any $c$.
 \end{proof}

 \begin{theorem} \label{Th:ChPoly}
Let $c = const.$, $n = mq + r$, $q, r \in \Integer$, $0 \le r \le m-1$.
Then the characteristic polynomial of the ideal $J_{n+1,c}$ is equal to
 \begin{equation}\label{Eq:ChP}
p_{2n,c} =  z^{(n)}\bar z^{(n)} + \sum_{k=1}^q (-1)^{mk}
\frac{c(c+1)\dots(c+k-1)}{(1-c)(2-c)\cdots (k-c)}\, z^{(n-mk;n+mk)}
 \end{equation}
where $z^{(p;q)} \bydef z^{(p)} \bar z^{(q)}+z^{(q)} \bar z^{(p)}$.
 \end{theorem}

 \begin{proof}
Use induction in $q$; base $q=0$. The ideal $J_{n+1,c}$ is
$I_2(m)$-invariant, and therefore $p_{2n,c}$ is $I_2(m)$-invariant, too.
So, $p_{2n,c}$ can be expressed as a polynomial of elementary invariants
$e_2 = z\bar z$ and $e_m = z^m + \bar z^m$. If $q = 0$, then the degree $2n
< 2m$, and $e_m$ cannot enter the expression. Therefore $p_{2n,c} = e_2^n$,
up to a multiple.

Induction step: let now $q > 0$. Normalize the polynomial $p_{2n,c}$ as in
\eqref{Eq:QuasiRes}. (Actually, this means multiplying all the terms of
\eqref{Eq:ChP} by $(1-c) \dots (q-c)$.) Then it is easy to see that if $r
\ne 0$, then the function $a_n(c)$ of Proposition \ref{Pr:Laplace} is equal
to $1$, while for $r = 0$ it is a polynomial of degree $1$ with leading
coefficient $1$, that is, $a_n(c) = c-\alpha$. So, for $c = \alpha$ one has
$\Delta p_{2n,c} = 0$, and therefore $z \bar z u \in J_{n+1,c}$ for any
polynomial $u \in \Complex_{2(n-1)}[z,\bar z]$. On the other hand, if $c =
\alpha$ is a regular value, then Lemma \ref{Lm:MultZZBar} implies that $u
\in J_{n,c}$ --- therefore, $\Complex_{2n-2}[z,\bar z] \subset J_{n,c}$,
which is impossible because $J_{n,c}$ is generated by two polynomials of
degree $n$. So, $\alpha$ is a singular value. Induction hypothesis and
Proposition \ref{Pr:Monom} show that $\alpha = q$ is the only possibility.

Thus, $a_n(c) = c-q$. This proves \eqref{Eq:ChP} by induction up to a term
$P$ annihilated by the operator $\Delta$. If $r \ne 0$, then this term is
zero because of $I_2(m)$-invariance. If $r = 0$, then it looks like $P =
\mu_n(c) (z^{(2n)} + \bar z^{(2n)})$ where $\mu_n$ is a polynomial of
degree $q$. Corollary \ref{Cr:Zeros} shows that $a_n(c)$ is divisible by
$(c+1) \dots (c+q-1)$.

Let now $c = 0$. In this case $F = \Delta$ by Proposition \ref{Pr:FExpl}, and
therefore $J_{n+1,c}$ is invariant under the action of the orthogonal group
$\name{SO}(2) = \{\xi \mid \lmod \xi\rmod = 1\} \subset \Complex^*$. So,
$p_{2n,0}$ is $\name{SO}(2)$-invariant, too, and therefore $p_{2n,0} = z^n
\bar z^n$ (as usual, up to a multiplicative constant). So, $\mu_n(0) = 0$,
and $\mu_n(c) = M_n c(c+1) \dots (c+q-1)$. Equation \eqref{Eq:QuasiRes}
allows to determine the normalizing coefficient: $M_n = 1$.
 \end{proof}

For $c \ne const.$ (and $m$ even) an explicit formula for $p_{2n,c}$ is
not yet known.

\section{Appendix: standard Frobenius algebras}\label{Sec:Apolar}

In this section we collect definitions and theorems about standard
Frobenius algebras used elsewhere in the paper. Our main sources are the
book \cite{Eisenbud} and the review \cite{Dol}.

\subsection{Main properties and operations}
 \begin{definition} \label{Df:StrFrob}
A graded commutative finite-dimensional $\Complex$-algebra $A =
\bigoplus_{k=0}^N A_k$ with a unit $\mathbf{1} \in A_0$ is called {\em
standard Frobenius} if
 \begin{enumerate}
\item\label{It:Socle} $\dim A_0 = 1$.
\item\label{It:Gener} $A$ is generated by the component $A_1$,
\item\label{It:Frob} The operation $(a,b) \mapsto \name{pr}_N(ab)$ defines
a non-degenerate bilinear form on $A$ (here $\name{pr}_N$ means taking the
homogeneous component of degree $N$).
 \end{enumerate}
 \end{definition}

A graded commutative algebra possessing properties \ref{It:Socle} and
\ref{It:Frob} only is called Frobenius; the bilinear form mentioned in
\ref{It:Frob} is called a Frobenius form. Algebras possessing property
\ref{It:Gener} are called standard. The number $N$ is called a degree of
the algebra $A$; the number $r = \dim A_1$ is called a rank of $A$.

Non-degeneracy of the Frobenius form implies symmetry of the Hilbert series
of a standard Frobenius algebra: $\dim A_k= \dim A_{N-k}$ for any $k =
0,1,\dots,N$; in particular, $\dim A_N = 1$.

For each standard Frobenius algebra $A=\bigoplus_{k=0}^N A_k$ one can
define its {\em characteristic polynomial} $p_A \in \Complex[A_1]$, which
appears to encode all the information about the algebra. It follows from
Property \ref{It:Gener} that $A \cong S(A_1)/I$ where $I = \bigoplus_{k\ge
0} I_k$ is a graded ideal. Since $\dim A_N = 1$, the component $I_N$ is a
co-dimension $1$ subspace of $S^N(A_1)$. Therefore it is a kernel of a
linear function $p_A \in (S^N(A_1))^*$. A natural identification
$(S^N(A_1))^* = S^N(A_1^*) = \Complex_N[A_1]$ allows to consider $p_A$ a
polynomial function of degree $N$ on $A_1$.

Clearly, $p_A$ is defined up to a scalar multiple; so, strictly speaking,
$p_A$ is an element of the projective space $\PP(\Complex_N[A_1])$.

 \begin{proposition} \label{Pr:DerivDep}
For every $k = 0, \dots, N$ the homogeneous component $A_k$ of the algebra
$A$ is equal to $S^N(A_1)/V_k$ where the subspace $V_k \subset S^k(A_1)$ is
the common kernel of all the partial derivatives of order $N-k$ of the
polynomial $p_A$.
 \end{proposition}

 \begin{proof}
Let $P \in S^k(A_1)$. Since the Frobenius form on $A$ is non-degenerate, one
has $P \in J_k$ if and only if $PQ \in J_N$ for all $Q \in S^{N-k}(A_1)$.
By the definition of $p_A \in (S^k(A_1))^*$, this is equivalent to $\langle
p_A, PQ\rangle = 0$. Fix a basis $x_1, \dots, x_r \in A_1$. Then $P \in
J_k$ if and only if for all $\alpha_1, \dots, \alpha_r \in \Integer_{\ge
0}$ such that $\alpha_1 + \dots + \alpha_r = N-k$ one has $0 = \langle p_A,
x_1^{\alpha_1} \dots x_r^{\alpha_r} P\rangle = \langle \frac{\partial^{N-k}
p_A}{\partial x_1^{\alpha_1} \dots \partial x_r^{\alpha_r}}, P\rangle$.
 \end{proof}

This result can be reformulated in a more elegant
coordinate-free manner. To do this, consider the action of the symmetric
algebra $S(V)$ on the algebra of polynomial functions $\Complex[V]$ by the
differential operators with constant coefficients. This action is graded:
$S^m(V) \times \Complex_n[V] \to \Complex_{n-m}[V]$ for any $n \ge m$.

 \begin{proposition}
Let $A = \bigoplus_{k=0}^N A_k$ be a standard Frobenius algebra with the
characteristic polynomial $p = p_A$. Then the dual space $A^* \subset
\Complex[A_1]$ is a cyclic $S(A_1)$-module $S(A_1)(p)$. More precisely,
$A_k^*=S^k(V)(p)$ for all $k = 0, \dots, N$.
 \end{proposition}

 \begin{corollary} \label{Cr:EssDep}
For any $v \in A_1$ the directional derivative $\pder{p_A}{v}$ does not
vanish identically. In other words, the partial derivatives
$\pder{p_A}{x_1}, \dots, \pder{p_A}{x_r}$ are linearly independent.
 \end{corollary}

This follows from the equality $\dim A_1 = r$. In still other words, the
polynomial $p_A$ depends essentially on $r$ variables and cannot be written
as a polynomial of a smaller number of variables. We denote by
$\Complex^0_N[V]$ the set of all $p\in \Complex_N[V]$ having this property.
Clearly, $\Complex^0_N[V]$ is an open dense $GL(V)$-invariant subset of
$\Complex_N[V]$.

 \begin{corollary}
As a graded vector space the algebra $A$ is isomorphic to the set of all
partial derivatives of the polynomial $p_A$.
 \end{corollary}

 \begin{theorem}
A graded isomorphism $f: A \to B$ of standard Frobenius algebras carries
$p_A$ to $p_B$. For any polynomial $p \in \Complex^0_N[V]$ there exists a
unique standard Frobenius algebra $A$ of degree $N$ such that $A_1 = V$ and
$p_A = p$.
 \end{theorem}

 \begin{proof}
The first assertion follows immediately from the definition of $p_A$. The
uniqueness in the second assertion is a consequence of Proposition
\ref{Pr:DerivDep}.

We now prove the existence of a standard Frobenius algebra with a given
$p_A = p$. Call an algebra {\em standard almost Frobenius} 
if it possesses properties \ref{It:Socle} and \ref{It:Gener} from 
Definition \ref{Df:StrFrob}. Note that for each standard quasi-Frobenius 
algebra $\tilde A$ the radical $\name{rad}(\tilde A) = \{v \in \tilde A 
\mid (v,w) = 0\,\forall w \in \tilde A\}$ of its Frobenius form 
$(\cdot,\cdot)$ is a graded ideal. One can easily show that 
$\name{rad}(\tilde A)$ is the maximal graded ideal in $\tilde A$ not 
intersecting the highest degree component $\tilde A_N$.

Now let $p \in \Complex_N^0[V]$. Let $I_N \subset S^N(V)$ be the kernel of
the corresponding linear function $S^N(V) \to \Complex$. Denote $J \bydef
I_N \oplus_{k>N} S^k(V)$. Clearly, this is an ideal in $S(V)$. Denote
$\tilde A \bydef S(V)/J$ and define a bilinear form on this algebra by the
formula $(v,w)_{\tilde A} = \langle p, uv\rangle$.

Clearly, $\tilde A$ is a standard quasi-Frobenius algebra. We define the
algebra $A$ as the quotient: $A \bydef \tilde A/\name{rad}(\tilde A)$. An
easy check shows that $A$ is standard Frobenius and $p_A = p$. Since $p \in
\Complex^0_N[V]$, the ideal $\name{rad}(\tilde A)$ does not intersect
$A_1$, and the equality $A_1 = V$ follows.
 \end{proof}

 \begin{remark}
This result implies that standard Frobenius algebras are in 1-1
correspondence with their characteristic polynomials modulo linear changes
of variables. In other words, the isomorphism classes of degree $N$
standard Frobenius algebras generated by $A_1 = V$ are parametrized by the
points of the orbit space $\PP\Complex_N^0[V]/GL(V)$ where
$\PP\Complex_N^0[V]$ is a dense subset of $\PP(\Complex_N[V])$
corresponding to $\Complex_N^0[V]$.
 \end{remark}

One can define several operations over the standard Frobenius algebras. One
of them is the tensor product:

 \begin{proposition}
The tensor product of standard Frobenius algebras $A$ and $B$ has a natural
structure of a standard Frobenius algebra with the characteristic
polynomial $p_{A\otimes B} = p_A\otimes p_B$.
 \end{proposition}

 \begin{proof}
If $A = \bigoplus_{k=0}^N A_k$, $B=\bigoplus_{l=0}^M B_k$, then the tensor
product $C = A \otimes B$ is a graded algebra: $C = \bigoplus_{r=0}^{M+N}
C_r$, where
 \begin{equation*}
C_r = \bigoplus A_k \otimes B_{r-k}
 \end{equation*}
The algebra $C$ is generated $C_1 = A_1 \otimes B_0 \oplus A_0 \otimes
B_1$. Also,
 \begin{equation*}
C_{M+N} = (C_1)^{M+N} = (A_1\otimes B_0 \oplus A_0\otimes B_1)^{M+N} =
A_N\otimes  B_M.
 \end{equation*}
It remains to prove that the bilinear form $(\cdot,\cdot)_C$ is
non-degenerate. Indeed, for any $a,a'\in A$, $b,b'\in B$:
 \begin{eqnarray*}
(a\otimes b,a'\otimes b') &=& \name{pr}_{C_{M+N}} (a\otimes b)(a' \otimes
b') = \name{pr}_{C_{M+N}}(aa' \otimes bb') \\
&=& \name{pr}_{A_N \otimes  B_M}(aa' \otimes bb') = \name{pr}_{A_N}(aa')
\name{pr}_{B_M}(bb') = (a,a')(b,b'),
 \end{eqnarray*}
that is, $(\cdot,\cdot)_{A \otimes B} = (\cdot,\cdot)_A \otimes (\cdot,\cdot)_B$, 
which implies non-degeneracy and the equality for characteristic
polynomials.
 \end{proof}

Another important operation is the internal product:

 \begin{definition}
Let $A = \bigoplus_{k=0}^N A_k$ be a standard Frobenius algebra and $U
\subset A_1$ be a subspace such that $U^N = A_N$. Define the {\em
restriction} $\left. A\right|_U$ of $A$ to $U$ to be the quotient $\langle
U\rangle/\mathcal{R}$ where $\langle U\rangle \subset A$ is the subalgebra
generated by $U$, and $\mathcal R$ is the radical of the restriction of the
Frobenius form $(\cdot,\cdot)$ to $\langle U\rangle$.
 \end{definition}

 \begin{definition}
Let $A$ and $B$ be standard Frobenius algebras such that $A_1=B_1=V$.
Define the {\em internal product} $A*B$ of the algebras as the restriction
$\left. A\otimes B\right|_{\Delta(V)}$ where $\Delta(V) = \{v\otimes
1+1\otimes v \mid v\in V\} \subset (A \otimes B)_1$ is the diagonal.
 \end{definition}

The following proposition is obvious:
 \begin{proposition}
Characteristic polynomials of standard Frobenius algebras satisfy the
identities $p_{\left.A\right|_U} = \left.p_A\right|_U$ and $p_{A*B} =
p_Ap_B$.
 \end{proposition}

 \begin{corollary}
The internal product is a monoidal operation on the category of standard
Frobenius algebras generated by a given space $V$.
 \end{corollary}

In other words, $A*(B*C) = (A*B)*C$ and $A*E = E*A = A$, where $E$ is a
trivial one-dimensional algebra.

\subsection{Free standard Frobenius algebras}

Recall from \cite{Eisenbud} that a graded zero-dimensional complete
intersection cone in $\Complex^\ell$ is a finite-dimensional quotient of 
$\Complex[x_1, \dots, x_\ell]$ by an an ideal $J$ generated by $\ell$ homogeneous polynomials $P_1, \dots,
P_\ell$.

 \begin{proposition}[\cite{Eisenbud}]\label{Pr:AllMinim2}
Each graded zero-dimensional complete intersection cone in $\Complex^\ell$
is a standard Frobenius algebra of degree $\deg A = \sum_{i=1}^\ell (\deg
P_i-1)$ and total dimension $\dim_\Complex A = \prod_{i=1}^\ell \deg P_i$.
 \end{proposition}

For the sake of simplicity, based on Proposition \ref{Pr:AllMinim2}, we
will refer to graded zero-dimensional complete intersection cones as {\em
free standard Frobenius algebras}. Borrowing a term from the Invariant
Theory, call the numbers $m_i = \name{deg}(P_i) - 1$ {\em exponents} of the
algebra $A$. It is easy to see that for each free standard Frobenius
algebra $A$ the elements $P_1, \dots, P_\ell$ are algebraically
independent. The Hilbert polynomial $h_A(t) = \sum_{k=0}^N \dim A_k t^k$ of
$A$ is
 \begin{equation}\label{eq:hilbert factor}
h_A(t) = \prod_{i=1}^\ell \frac{1-t^{m_i+1}}{1-t}.
 \end{equation}

Here are some examples of free standard Frobenius algebras.

 \begin{example}
We say that a standard Frobenius algebra $A$ is monomial if there exists a
basis $e_1, \dots, e_\ell$ in $A_1$ such that $p_A = e_1^{m_1} e_2^{m_2}
\dots e_\ell^{m_\ell}$ for some integers $1 \le m_1 \le m_2 \le \dots \le
m_\ell$ such that $m_1 + m_2 + \dots + m_\ell = N$.

Each monomial standard Frobenius algebra $A$ with the characteristic
polynomial $p_A$ of this form is free with the exponents $m_1, m_2, \dots,
m_\ell$. Moreover, $A\cong \bigotimes_{i=1}^\ell A(m_i)$ where $A(m) \bydef
\Complex[x]/(x^{m+1})$.
 \end{example}

 \begin{example}
Another example of a free standard Frobenius algebra is the algebra of
coinvariants of a reflection group. Let $W$ be a finite reflection group of
a space $V$. Denote by $I_W$ the ideal in $S(V)$ generated by all
homogeneous $W$-invariant elements of positive degrees. It is well-known
(see \cite{Dol}) that the algebra $A_W = S(V)/I_W$ of coinvariants is
standard Frobenius.

Let $\Delta = \Delta_W\in \Complex[V]$ be the polynomial  skew-invariant of $W$ of the smallest degree. Clearly, $\Delta$ is unique up to a constant multiple and factors into the product of all positive roots of $W$. If $W$ is irreducible, then the degree of $\Delta$ is
equal to $\lmod \Refl\rmod = h\ell/2$, where $h$ is the Coxeter number of
$W$ and $\Refl$ is the set of all reflections in $W$ (see Section
\ref{SSec:Defs}).

 \begin{proposition}[\cite{Dol}]
The coinvariant algebra $A_W$ is free standard Frobenius with the
characteristic polynomial $\Delta$.
 \end{proposition}
 \end{example}

 \begin{example}
Finite-dimensional quotients of the module $M(\triv)$ over a rational
Cherednik algebra $H_c(W)$ (see Section \ref{SSec:Defs} above for notation)
are free standard Frobenius algebras. According to \cite{BEG}, if one
identifies $M(\triv) = S(V)$ and considers the module in question as a
quotient algebra, then this algebra is free standard Frobenius. See Section
\ref{SSec:Mult} for a complete description of this algebra for $W =
I_2(m)$.
 \end{example}

\subsection{Rank $2$ case}
We  refine Proposition \ref{Pr:AllMinim2} when $\dim
A_1=2$ as follows.

\begin{proposition}[\cite{Eisenbud}]
Let $A$ be a quotient algebra of $\Complex[x_1,x_2]$ by a homogeneous ideal
$J$. Then the following are equivalent:
 \begin{enumerate}
\item $A$ is standard Frobenius.
\item $A$ is free standard Frobenius.
\item The ideal $J$ is generated by two coprime homogeneous polynomials
$R_1, R_2$.
 \end{enumerate}
\end{proposition}

For such algebras it is possible to write down an explicit expression for
the characteristic polynomial provided $R_1, R_2$ are known.

Let $A = \Complex[x_1, x_2]/J$ be a standard Frobenius algebra, and the
ideal $J$ be generated by coprime polynomials $R_1, R_2 \in
\Complex_{n+1}[x_1, x_2]$. Then by Proposition \ref{Pr:AllMinim2} the degree
of $A$ is $N = 2n$, and the component $J_N$ of the ideal $J$ is spanned by
all the polynomials $x_1^\alpha x_2^\beta R_i$ where $\alpha+\beta = n-1$
and $i = 1,2$. If $R_1 = \sum_{i=0}^{n+1} a_i x_1^i x_2^{n+1-i}$ and $R_2 =
\sum_{i=0}^{n+1} b_i x_1^i x_2^{n+1-i}$, then the functional $p_A: \Complex_{2n}[x_1,x_2] \to \Complex$ should
vanish on all the polynomials $\sum_{i=0}^{n+1} a_i x_1^{i+\alpha}
x_2^{2n-i-\alpha}$ and $\sum_{i=0}^{n+1} b_i x_1^{i+\alpha}
x_2^{2n-i-\alpha}$ where $\alpha = 0, \dots, n-1$. It means that
 \begin{equation}\label{Eq:QuasiRes}
p_A = \sum_{i=0}^{2n} u_i y_1^{(i)} y_2^{(2n-i)}
 \end{equation}
where $a^{(b)} \bydef a^b/b!$, and $u_i$ is the minor of the $(2n+1)\times
(2n)$ matrix $U$ composed of the coefficients of these polynomials:
 \begin{equation*}
U = \left(\begin{array}{cccccccc}
a_0 & a_1 & \dots & a_n & a_{n+1}     & 0   & \dots & 0\\
0   & a_0 & \dots & a_{n-1} & a_n & a_{n+1} & \dots & 0\\
\vdots  \\
0   & 0   & \dots & a_0     & a_1     & a_2 & \dots & a_{n+1}\\
b_0 & b_1 & \dots & b_n & b_{n+1}     & 0   & \dots & 0\\
0   & b_0 & \dots & b_{n-1} & b_n & b_{n+1} & \dots & 0\\
\vdots  \\
0   & 0   & \dots & b_0     & b_1     & b_2 & \dots & b_{n+1}
 \end{array}\right)
 \end{equation*}
obtained by deletion of its $(i+1)$-th column.

Let $p = \sum_{i=0}^n c_i x_1^{(i)} x_2^{(n-i)} $ be a homogeneous
polynomial (a binary form) of degree $N$. Define the matrix $H_k(p)$ by $(H_{k,n-k}(p))_{ij} = c_{i+j}$ 
for $i = 0,1,\dots,k$, $j = 0,1,\dots,n-k$ (see e.g.\ \cite{Elk}). Clearly,
the transposed of $H_{n-k,k}(p)$ is  $H_{k,n-k}(p)$ for all $k \le n$.

 \begin{proposition}
Let $A$ be a rank $2$ standard Frobenius algebra with the characteristic
polynomial $p$, and let $k \le \frac{n}{2}$. Then $\dim A_k = k+1$ if and
only if the matrix $H_{k,n-k}(p)$ has the maximal rank $k+1$.
 \end{proposition}

The proof follows from \eqref{Eq:QuasiRes} and Proposition \ref{Pr:DerivDep}.


\begin{thebibliography}{99}
 
\bibitem{BEG0} Y.~Berest, P.~Etingof, V.~Ginzburg, Cherednik algebras and 
differential operators on quasi-invariants, {\em Duke Math. J.} {\bf 118} 
(2003), pp.~279--337.

\bibitem{BEG} Y.~Berest, P.~Etingof, V.~Ginzburg, Finite-dimensional
representations of rational Cherednik algebras. {\em Int.\ Math.\ Res.\
Not.} (2003), no.~19, pp.~1053--1088.

\bibitem{BrMaRo} M.~Broue, G.~Malle, R.~Rouquier,  Complex reflection 
groups, braid groups, Hecke algebras, {\em J. Reine Angew. Math.} {\bf 500} 
(1998), 127-190.

\bibitem{Cherednik}  I.~Cherednik, Double affine Hecke algebras and
Macdonald's conjectures. {\em Ann.\ of Math.} (2) {\bf 141} (1995), no.~1,
pp.~191--216.

\bibitem{Chev} C.~Chevalley, Invariants of finite groups generated by
reflections, {\em Amer.\ J.\ Math.} {\bf 77} (1955), pp.~778--782.

\bibitem{Chmutova} T.~Chmutova, Representations of the rational
Cherednik algebras of dihedral type, {\em Journal of Algebra} 
{\bf 297}, No. 2, 2006, pp.~542--565.

\bibitem{CE} T.~Chmutova, P.~Etingof, On some representations of the 
rational Cherednik algebra,  {\em Represent. Theory} {\bf 7} (2003), 
641--650 (electronic).

\bibitem{Dez} Ch.~Dez\'el\'ee, Repr\'esentations de dimension finie de
l'alg\`ebre de Cherednik rationnelle. {\em Bull.\ Soc.\ Math.\ France} {\bf
131} (2003), no.~4, pp.~465--482.

\bibitem{Dol} I.~Dolgachev, {\em Lectures on invariant theory}. London
Mathematical Society Lecture Note Series, 296. Cambridge University Press,
Cambridge, 2003.

\bibitem{DunklJeuOpdam} C.~Dunkl, M.~de Jeu, E.~Opdam, Singular
polynomials for finite reflection groups. {\em Trans. Amer. Math. Soc}.
{\bf 346} (1994), no.~1, pp.~237--256.

\bibitem{DunklDef} Differential-difference operators associated to
reflection groups, {\em Trans.\ Amer.\ Math. Soc.} {\bf 311} (1989),
pp.~167--183.

\bibitem{DO} C.~Dunkl, E.~Opdam, Dunkl operators for complex reflection
groups, {\em Proc.\ London Math.\ Soc.} (3) {\bf 86} (2003), no.~1,
pp.~70--108.

\bibitem{DuInter} C.~Dunkl, Intertwining operators and polynomials
associated with symmetric group, {\em Monatshefte Math.} {\bf 126} (1998),
pp.~181--209.

\bibitem{DuSing} C.~Dunkl,  Singular polynomials for the symmetric
groups,  {\em Int. Math. Res. Not.} 2004, no. {\bf 67}, 3607--3635.

\bibitem{DuSingMod} C.~Dunkl, Singular polynomials and modules for
the symmetric groups, {\em Int. Math. Res. Not.}, 2005, {\bf 39}, pp.~2409--2436.

\bibitem{Eisenbud} D.~Eisenbud, {\em Commutative algebra with a view
toward algebraic geometry}, Graduate Texts in Mathematics, 150.
Springer-Verlag, NY, 1995.

\bibitem{Elk} N.~Elkies, On finite sequences satisfying linear recursions.
{\em New York J. Math.} {\bf 8} (2002), pp.~85--97 (electronic)

\bibitem{GorQuot} I.~Gordon On the quotient by diagonal invariants, {\em
Invent.\ Math.}, {\bf 153} (2003), pp.~503--518.

\bibitem{GorI} I.~Gordon, J.~Stafford, Rational Cherednik algebras
and Hilbert schemes I, {\em Adv. Math.}, {\bf 198} (2005), no. 1, pp.~222--274.

\bibitem{GorII} I.~Gordon, J.T.~Stafford, Rational Cherednik
algebras and Hilbert schemes II: representations and sheaves, {\em Duke Math. J.}, {\bf 132} (2006), no. 1, pp.~73-135.

\bibitem{Haiman} M.~Haiman, Conjectures on the quotient ring by diagonal
invariants, {\em J.\ Algebraic Combin.} {\bf 3} (1994), no. 1, pp.~17--76.

\bibitem{Hum} J.~Humphreys, {\em Reflection groups and Coxeter groups},
Cambridge Studies in Advanced Mathematics, {\bf 29}, Cambridge University
Press, Cambridge, 1990.

\bibitem{Kos} B.~Kostant, Clifford algebra analogue of the 
Hopf-Koszul-Samelson theorem, the $\rho$-decomposition ${\mathcal 
C}(\mathfrak{g}) = \name{End} V_\rho \otimes {\mathcal C}(P)$, and the 
$\mathfrak{g}$-module structure of $\bigwedge \mathfrak{g}$, {\em Adv. 
Math.}, {\bf 125} (1997), no. 2, pp.~275--350.

\bibitem{Rouq} R.~Rouquier, Representations of rational Cherednik algebras, 
in ``Infinite-dimensional aspects of representation theory and applications,''  {\em American Math. Soc.}, 2005, pp.~103-131.

\bibitem{SheTodd} G.~Shephard, J.~Todd, Finite unitary reflection
groups, {\em Canad.\ J.\ Math.} {\bf 6} (1954), pp.~274--304.

 \end{thebibliography}
\end{document}